\documentclass[mathpazo,online]{jml}

\renewcommand\thisvolume{4}
\renewcommand\thisnumber{4}
\renewcommand\thisyear {2025}
\renewcommand\thismonth{December}

\renewcommand\originalreceiveddate{13/02/2025}

\renewcommand\acceptdate{08/08/2025}

\renewcommand\pagestart{223}
\renewcommand\pageend{263}

\renewcommand\articledoi{doi.org/10.4208/jml.250213}

\renewcommand\handler{xxx}

\usepackage{lipsum}

\usepackage[utf8]{inputenc}
\usepackage{mathrsfs,amsthm,xcolor,verbatim,bbm,amsmath,amsfonts,
amssymb,nicefrac,enumitem}
\usepackage{hyperref,bm,mathtools,xparse,etoolbox}
\usepackage[capitalise,sort]{cleveref} %

\newcommand{\revision}[1]{#1}

\newcommand{\revisionjmla}[1]{#1}
\newcommand{\revisionjmlb}[1]{#1}

\usepackage{textcomp}
\usepackage{sidecap}

\DeclareMathOperator*{\argmin}{arg\,min}

\crefname{enumi}{item}{items}

\crefname{figure}{Figure}{Figures}
\crefname{subsection}{Subsection}{Subsections}

\crefname{case}{Case}{Cases}
\crefname{cor}{Corollary}{Corollaries}

\renewcommand{\emptyset}{\varnothing}

\usepackage{tikz}
\usetikzlibrary{matrix,chains,positioning,decorations.pathreplacing,arrows}
\usetikzlibrary{shapes,arrows}
\tikzset{
	font={\fontsize{9pt}{12}\selectfont}}

\ExplSyntaxOn

\bool_new:N \g_noteobserve

\NewDocumentCommand{\setnote}{}{
  \bool_gset_true:N \g_noteobserve
}

\NewDocumentCommand{\setobserve}{}{
  \bool_gset_false:N \g_noteobserve
}

\NewDocumentCommand{\nobs}{ o }{
  \IfValueT{#1}{
    \str_if_eq:noTF {note} {#1} {
      \bool_gset_true:N \g_noteobserve
    } {
      \str_if_eq:noTF {Note} {#1} {
        \bool_gset_true:N \g_noteobserve
      } {
        \bool_gset_false:N \g_noteobserve
      }
    }
  }
  \bool_if:nTF { \g_noteobserve } {
    \bool_gset_false:N \g_noteobserve
    note
  } {
    \bool_gset_true:N \g_noteobserve
    observe
  }
  \IfValueF{#1}{~}
}

\NewDocumentCommand{\Nobs}{ o }{
  \IfValueT{#1}{
    \str_if_eq:noTF {note} {#1} {
      \bool_gset_true:N \g_noteobserve
    } {
      \str_if_eq:noTF {Note} {#1} {
        \bool_gset_true:N \g_noteobserve
      } {
        \bool_gset_false:N \g_noteobserve
      }
    }
  }
  \bool_if:nTF { \g_noteobserve } {
    \bool_gset_false:N \g_noteobserve
    Note
  } {
    \bool_gset_true:N \g_noteobserve
    Observe
  }
  \IfValueF{#1}{~}
}

\int_new:N \g_furthermore

\NewDocumentCommand{\Moreover}{ o o }{
  \IfValueT{#1}{
    \str_case:nn {#1} {
      {Furthermore} {\int_set:Nn {\g_furthermore} {0}}
      {Moreover} {\int_set:Nn {\g_furthermore} {1}}
      {In~addition} {\int_set:Nn {\g_furthermore} {2}}
      {note} {\bool_gset_true:N \g_noteobserve}
      {observe} {\bool_gset_false:N \g_noteobserve}
    }
    \IfValueT{#2}{
      \str_case:nn {#2} {
        {Furthermore} {\int_set:Nn {\g_furthermore} {0}}
        {Moreover} {\int_set:Nn {\g_furthermore} {1}}
        {In~addition} {\int_set:Nn {\g_furthermore} {2}}
        {note} {\bool_gset_true:N \g_noteobserve}
        {observe} {\bool_gset_false:N \g_noteobserve}
      }
    }
  }
  \int_case:nn { \int_mod:nn {\g_furthermore} {3} } {
    { 0 } { Furthermore,~\nobs that}
    { 1 } { Moreover,~\nobs that}
    { 2 } { In~addition,~\nobs that}
  }
  \int_incr:N \g_furthermore
  \IfValueF{#1}{~}
}

\bool_new:N \g_hencetherefore

\NewDocumentCommand{\hence}{}{
  \bool_if:nTF { \g_hencetherefore } {
    \bool_gset_false:N \g_hencetherefore
    hence~
  } {
    \bool_gset_true:N \g_hencetherefore
    therefore~
  }
}

\NewDocumentCommand{\Hence}{}{
  \bool_if:nTF { \g_hencetherefore } {
    \bool_gset_false:N \g_hencetherefore
    Hence,~we~obtain~
  } {
    \bool_gset_true:N \g_hencetherefore
    Therefore,~we~obtain~
  }
}

\seq_new:N \g_cflist_loaded
\seq_new:N \g_cflist_pending

\NewDocumentCommand{\cfadd}{ m }
{
	\seq_if_in:NnF \g_cflist_loaded { #1 } {
		\seq_if_in:NnF \g_cflist_pending { #1 } {
			\seq_gput_right:Nn \g_cflist_pending { #1 }
		}
	}
}

\NewDocumentCommand{\cfconsiderloaded}{ m }{
	\seq_gput_right:Nn \g_cflist_loaded {#1}
}

\NewDocumentCommand{\cfremove}{ m }
{
	\seq_gremove_all:Nn \g_cflist_pending { #1 }
}

\NewDocumentCommand{\cfload}{ o }
{
	\seq_if_empty:NTF \g_cflist_pending {\unskip} {
		(cf.\ \cref{\seq_use:Nn \g_cflist_pending {,}})\IfValueTF{#1}{#1~}{\unskip}
		\seq_gconcat:NNN \g_cflist_loaded \g_cflist_loaded \g_cflist_pending
		\seq_gclear:N \g_cflist_pending
	}
}

\NewDocumentCommand{\cfclear} {} {
	\seq_gclear:N \g_cflist_loaded
	\seq_gclear:N \g_cflist_pending
}

\NewDocumentCommand{\cfout}{ o }
{
	\seq_if_empty:NTF \g_cflist_pending {\unskip} {
		(cf.\ \cref{\seq_use:Nn \g_cflist_pending {,}})\IfValueTF{#1}{#1~}{\unskip}
		\seq_gclear:N \g_cflist_pending
	}
}

\NewDocumentCommand{\ifnocf} { m } {
	\seq_if_empty:NT \g_cflist_pending { #1 }
}

\ExplSyntaxOff

\NewDocumentEnvironment{cproof}{m}
{\begin{proof}[Proof of \cref{#1}]}%
	{\noindent The proof of \cref{#1} is thus complete.
\end{proof}}

\NewDocumentEnvironment{cproof2}{m}
{\begin{proof}[Proof of \cref{#1}]}%
	{\noindent This completes the proof of \cref{#1}.
\end{proof}}

\newcommand{\citep}{\cite}
\newcommand{\citet}{\cite}
\usepackage{url}
\usepackage{enumitem}
\usepackage{mathtools}
\usepackage{adjustbox}
\usepackage{floatrow}
\newfloatcommand{capbtabbox}{table}[][\FBwidth]
\usepackage{blindtext}

\usepackage{amsmath,amsfonts,bm,mathtools,amsthm,amssymb,commath,cases}
\usepackage{algorithm}
\usepackage{algpseudocode} 
\usepackage{graphicx}
\usepackage{subcaption}
\usepackage{comment}
\usepackage{xspace}
\usepackage{booktabs}
\usepackage{caption}
\usepackage[normalem]{ulem}
\usepackage{enumitem}
\setlist[itemize]{topsep=3pt, itemsep=2pt}

\newcommand{\bR}{\mathbb{R}}

\newcommand{\Var}{\mathrm{Var}}

\usepackage{colortbl}

\usepackage{makecell}

\theoremstyle{plain}
\newtheorem{thm}{Theorem}

\makeatletter
\DeclareRobustCommand\onedot{\futurelet\@let@token\@onedot}
\def\@onedot{\ifx\@let@token.\else.\null\fi\xspace}

\def\eg{\emph{e.g}\onedot} 
\def\ie{\emph{i.e}\onedot}

\makeatother

\newcommand{\rmd}{\mathrm{d}}
\newcommand{\LQR}{\text{LQR}}

\newcommand{\QR}{\text{QR}}
\newcommand{\tmin}{\text{min}}
\newcommand{\tmax}{\text{max}}
\newcommand{\xinit}{x_{\textrm{init}}}
\newcommand{\txinit}{\tilde{x}_{\textrm{init}}}
\newcommand{\transpose}{^{\operatorname{T}}}

\xspaceaddexceptions{1 2 3}
\def\dagger{DAgger\xspace}
\def\idealRatio{with the ideal curve being a straight horizontal segment passing ratio = $1$, percentage=$100\%$}

\begin{document}

\title{Progressive Optimal Path Sampling for Closed-Loop Optimal Control Design with Deep Neural Networks
\thanks{An earlier version of this work appeared on arXiv under the title ``Initial Value Problem Enhanced Sampling for Closed-Loop Optimal Control Design with Deep Neural Networks
''.}
}

\author[1]{
Xuanxi Zhang
	\thanks{These authors contributed equally.
{\tt xuanxizhang@nyu.edu}.
	}
}
\author[2]{
Jihao Long
	\thanks{These authors contributed equally.
		{\tt longjh1998@gmail.com}.
	}
}
\author[2]{Wei Hu\thanks{{\tt weihu.math@gmail.com}}}
\author[3,4,5]{Weinan E \thanks{{\tt weinan@math.pku.edu.cn}}}
\author[6]{Jiequn Han \thanks{Corresponding author. jiequnhan@gmail.com}}

\affil[1]{Courant Institute of Mathematical Sciences, New York University}
\affil[2]{Institute for Advanced Algorithms Research, Shangha}
\affil[3]{AI for Science Institute, 100080, Beijing, P.R. China}
\affil[4]{School of Mathematical Science, Peking University, 100871, Beijing, P.R. China}
\affil[5]{Center for Machine Learning Research, Peking University, 100871, Beijing, P.R. China}
\affil[6]{Center for Computational Mathematics, Flatiron Institute}

\begin{abstract}

Closed-loop optimal control design for high-dimensional nonlinear systems has been a long-standing challenge. Traditional methods, such as solving the associated Hamilton-Jacobi-Bellman equation, suffer from the curse of dimensionality. Recent literature proposed a new promising approach based on supervised learning, by leveraging powerful open-loop optimal control solvers to generate training data and neural networks as efficient high-dimensional function approximators to fit the closed-loop optimal control. This approach successfully handles certain high-dimensional optimal control problems but still performs poorly on more challenging problems.  One of the crucial reasons for the failure is the so-called distribution mismatch phenomenon brought by the controlled dynamics. In this paper, we investigate this phenomenon and propose the Progressive Optimal Path Sampling (POPS) method to mitigate this problem. We theoretically prove that this enhanced sampling strategy outperforms both the vanilla approach and the widely used Dataset Aggregation (DAgger) method on the classical linear-quadratic regulator by a factor proportional to the total time duration. We further numerically demonstrate that the proposed sampling strategy significantly improves the performance on tested control problems, including the optimal landing problem of a quadrotor and the optimal reaching problem of a 7 DoF manipulator.
\end{abstract}

\keywordone{Closed-loop optimal control,}
\keywordtwo{Distribution mismatch, }
\keywordthree{Adaptive sampling, }
\keywordfour{Supervised learning.}

\maketitle

\section{Introduction}
\revisionjmla{Optimal control aims to find a control for a dynamical system over a period of time such that a specified cost function is minimized. This cost often reflects a combination of task-specific goals such as energy usage, deviation from a tracking target, and control effort. Finding such an optimal control should be distinguished from classical stabilization control~\citep{franklin2002feedback}, which focuses only on keeping the system state bounded or driving it to an equilibrium, without regard to minimizing a cost.}
Generally speaking, there are two types of optimal controls: open-loop optimal control and closed-loop (feedback) optimal control. Open-loop optimal control, also known as trajectory optimization, deals with the problem with a given initial state,
and its solution is a function of time for the specific initial data, independent of the other states of the system. 
In contrast, closed-loop optimal control aims to 
find the optimal control policy as a function  of  the state that gives us optimal control for general initial states.

By the nature of the problem, solving the open-loop control problem is relatively easy and various open-loop control solvers can handle nonlinear problems even when the state lives in high dimensions \citep{betts1998survey,rao2009survey}. Closed-loop control is much more powerful than open-loop control since it can cope with different initial states, and it is more robust to the disturbance of dynamics. The classical approach to obtaining a closed-loop optimal control function is by solving the associated Hamilton-Jacobi-Bellman (HJB) equation. However, traditional numerical algorithms for HJB equations such as the finite difference method or finite element method face the curse of dimensionality \citep{bellman1957} and hence can not deal with high-dimensional problems.

\revision{There is a long history of employing neural networks (NNs) to solve the optimal control problems; see \eg  \citet{NNHJB2005,NNValue2007,GPS2013,VariationalPolicySearchTO2013,levine2014learning,combiningTO2014,han2016deep,PLATO2016,APEX2021,nakamura2021adaptive,nakamura2021neural,qrnet1,bottcher2022aidirect,e2022empowering}, and it is getting more attention recently since neural networks have demonstrate superior representation and generalization capabilities.}

\revision{Generally speaking, there are two categories of methods in this promising direction. One is direct policy search approach~\citep{GPS2013,combiningTO2014,han2016deep,pmlr-v144-ainsworth21a,bottcher2022aidirect,pmlzhao22advreg},  
which parameterizes the policy function by NNs, samples the total cost with various initial points, and directly minimizes the average total cost.
When learning complex policies with hundreds of parameters or solving problems with a long time span and high nonlinearity, the corresponding optimization problems can be extremely hard and may get stuck in poor local minima \citep{GPS2013,zhao2022offline}.
The other category of methods is based on supervised learning \citep{levine2014learning,combiningTO2014,nakamura2021adaptive,nakamura2021neural,qrnet1,nakamura2022neural}. 
Combining various techniques for open-loop control, one can solve complex high-dimensional open-loop optimal control problems; see \citet{betts1998survey,rao2009survey,kang2021algorithms} for surveys.
Consequently, we can collect optimal trajectories for different initial points as training data, parameterize the control function (or value function) using NNs, and train the NN models to fit the closed-loop optimal controls (or optimal values). This work focuses on the second approach and aims to improve its performance through adaptive sampling.}

As demonstrated in \citet{nakamura2021neural,zang2022machine,zhao2022offline}, NN controllers trained by the vanilla supervised-learning-based approach can perform poorly even when both the training error and test error on collected datasets are fairly small. Some existing works attribute this phenomenon to the fact that the learned controller may deteriorate badly at some difficult initial states even though the error is small in the average sense. Several adaptive sampling methods regarding the initial points are hence proposed (see Section \ref{sec:compare} for a detailed discussion). 
However, these methods all focus on choosing optimal paths according to different initial points and ignore the effect of dynamics.  This is an issue since the paths controlled by the NN will deviate from the optimal paths further and further over time due to the accumulation of errors. 
As shown in Section \ref{sec:prob_quadrotor}, applying adaptive sampling only on initial points is insufficient to solve challenging problems.

This work focuses on the so-called \textit{distribution mismatch} or \textit{covariance shift} phenomenon brought by the dynamics in the supervised-learning-based approach. This phenomenon refers to the fact that the discrepancy between the state distribution of the training data and the state distribution generated by the NN controller typically increases over time and the training data fails to represent the states encountered when the trained NN controller is used. Such phenomenon has also been identified in reinforcement learning \citep{kakade2002approximately,long2022perturbational} and imitation learning \citep{ross2010efficient,ross2011reduction}.

\revision{In this paper, we propose the \textbf{Progressive Optimal Path Sampling (POPS)} method to update the states in the training dataset to more closely match the states the controller reaches.
While similar in spirit to the dataset aggregation (\dagger) method~\cite{ross2011reduction} in imitation learning, POPS addresses issues of applying \dagger to the optimal control problem considered, particularly the costs of solving open-loop control problems.
Given the substantial computational demands of querying the expert policy, it is not feasible to employ the mixed policies strategy of \dagger. Instead, our approach involves progressively sampling new states by the current neural network, deliberately avoiding long-horizon rollouts that exceed the network's capability.
Furthermore, unlike approaches that use model predictive control (MPC) as expert policy~\citep{PLATO2016,APEX2021}, which utilize only the first value from a trajectory, our method incorporates entire trajectories into the training dataset, reflecting the strength of solving open-loop control problems to produce a full trajectory of optimal labels.}

The resulting supervised-learning-based approach empowered by POPS can be interpreted as an instance of the exploration-labeling-training (ELT) algorithms \citep{zhang2018reinforced,e2021mlmodeling} for closed-loop optimal control problems.
At a high level, the ELT algorithm proceeds iteratively with the following three steps: (1) exploring the state space and examining which states need to be labeled; (2) solving the control problem to label these states and adding them to the training data; (3) training the machine learning model.
Through the lens of the ELT algorithm, there are at least three aspects to improve the efficiency of the supervised-learning-based approach for the closed-loop optimal control problem:

\begin{itemize}
\vspace{-0.2cm}
\item Use the adaptive sampling method. Adaptive sampling methods aim to sequentially choose the time-state pairs based on previous results to improve the performance of the NN controller. \textit{This corresponds to the first step in the ELT algorithm and is the main focus of this work.} We will discuss other adaptive sampling methods in Section \ref{sec:compare}.

\item Improve the efficiency of data generation, \ie, solving the open-loop optimal control problems. Although the open-loop optimal control problem is much easier than the closed-loop optimal control problem, its time cost cannot be neglected and the efficiency varies significantly with different methods. This corresponds to the second step in the ELT algorithm and we refer to \citet{kang2021algorithms} for a detailed survey.

\item Improve the learning of the machine learning model. This corresponds to the third step in the ELT algorithm. The recent works \citet{qrnet1,nakamura2021neural,nakamura2022neural,hu2023learning} design a special network architecture such that the NN controller is close to the linear quadratic controller around the equilibrium point to improve the stability of the NN controller. Note that besides the popular use of neural networks, other models have been explored for approximating control/value functions~\citep{atkeson2002nonparametric,atkeson2007random,coates2008learning,deisenroth2009gaussian}, and the adaptive sampling methods developed in this paper are also applicable to these models.
\end{itemize}
\vspace{-0.2cm}

The main contributions of the paper can be summarized as follows:
\begin{itemize}
    \item  We investigate the distribution mismatch phenomenon brought by the controlled dynamics in the supervised-learning-based approach, which explains the failure of this approach for challenging problems. We propose POPS as an enhanced sampling method to update the training data, which significantly alleviates the distribution mismatch problem.
\item We show that POPS can significantly improve the performance of the learned closed-loop controller on a uni-dimensional linear quadratic control problem (theoretically and numerically) and two high-dimensional problems (numerically), the quadrotor landing problem and the reaching problem of a 7-DoF manipulator. 
\item We compare POPS with other adaptive sampling methods and show that POPS gives the best performance.
\end{itemize}

\vspace{-0.2cm}

\section{Preliminary}
\subsection{Open-loop and Closed-loop Optimal Control}
We consider the following deterministic controlled dynamical system:
\begin{equation}\label{eq:dyanmic_system}
\begin{dcases}
\dot{\bm{x}}(t) = \bm{f}(t,\bm{x}(t),\bm{u}(t)),\;t \in [t_0,T],\\
\bm{x}(t_0) = \bm{x}_0,
\end{dcases}
\end{equation}

where  $\bm{x}(t) \in \bR^n$ denotes the state, $\bm{u}(t) \in \mathcal{U} \subset \bR^m$ denotes the control with $\mathcal{U}$ being the set of admissible controls, $\bm{f}\colon [0,T]\times \bR^n\times \mathcal{U} \rightarrow \bR^n$ is a smooth function describing the dynamics, $t_0 \in [0,T]$ denotes the initial time, and $\bm{x}_0 \in \bR^n$ denotes the initial state. Given a fixed $t_0 \in [0,T]$ and $\bm{x}_0 \in \bR^n$, solving the open-loop optimal control problem means to find a control path $\bm{u}^*\colon  [t_0,T]\rightarrow \mathcal{U}$ to minimize
\begin{equation}\label{eq:loss}
 J(\bm{u};t_0,\bm{x}_0) =  \int_{t_0}^T L(t,\bm{x}(t),\bm{u}(t))\rmd t + M(\bm{x}(T))\,\,\text{ s.t. }(\bm{x},\bm{u}) \text{ satisfy the system \eqref{eq:dyanmic_system}}, \end{equation}
where $L\colon [0,T]\times\bR^n\times \mathcal{U}\rightarrow \bR$ and $M\colon \bR^n \rightarrow \bR$ are the running cost and terminal cost, respectively. We use $\bm{x}^*(t;t_0,\bm{x}_0)$ and $\bm{u}^*(t;t_0,\bm{x}_0)$ to denote the optimal state and control with the specified initial time $t_0$ and initial state $\bm{x}_0$, which emphasizes the dependence of the open-loop optimal solutions on the initial time and state. We assume the open-loop optimal control problem is well-posed, \ie, the solution always exists and is unique.

In contrast to the open-loop control being a function of time only, closed-loop control is a function of the time-state pair $(t,\bm{x})$. Given a closed-loop control $\bm{u}\colon [0,T]\times \mathbb{R}^n \rightarrow \mathcal{U}$, we can induce a family of the open-loop controls with all possible initial time-state pairs $(t_0,\bm{x}_0)$:
$
    \bm{u}(t;t_0,\bm{x}_0) = \bm{u}(t,\bm{x}_{\bm u}(t;t_0,\bm{x}_0)),
$
where $\bm{x}_{\bm u}(t;t_0,\bm{x}_0)$ is defined by the following \textit{initial value problem} (IVP):
\begin{equation}\label{eq:IVP}
\text{IVP}(\bm x_0, t_0, T, \bm{u}):
\begin{dcases}
\dot{\bm{x}}_{\bm u}(t;t_0,\bm{x}_0) = \bm{f}(t,\bm{x}_{\bm u}(t;t_0,\bm{x}_0),\bm{u}(t,\bm{x}_{\bm u}(t;t_0, \bm{x}_0)),\; t \in [t_0,T],\\
\bm{x}_{\bm u}(t_0;t_0,\bm{x}_0) = \bm{x}_0.
\end{dcases}
\end{equation}
To ease the notation, we always use the same character to denote the closed-loop control function and the induced family of the open-loop controls. The context of closed-loop or open-loop control can be inferred from the arguments and will not be confusing. It is well known in the classical optimal control theory (see, \eg \citet{liberzon2011calculus}) that there exists a closed-loop optimal control function $\bm{u}^*\colon [0,T] \times \mathbb{R}^n \rightarrow \mathcal{U}$ such that for any $t_0 \in [0,T]$ and $\bm{x}_0 \in \bR^n$,
$
\bm{u}^*(t;t_0,\bm{x}_0) = \bm{u}^*(t,\bm{x}^*(t;t_0,\bm{x}_0)),
$
which means the family of the open-loop optimal controls with all possible initial time-state pairs can be induced from the closed-loop optimal control function.
Since IVPs can be easily solved, one can handle the open-loop control problems with all possible initial time-state pairs if a good closed-loop control solution is available.
Moreover, the closed-loop control is more robust to dynamic disturbance and model misspecification, and hence it is much more powerful in applications. In this paper, our goal is to find a near-optimal closed-loop control $\hat{\bm{u}}$ such that  for $\bm{x}_0 \in X \subset \bR^n$ with $X$ being the set of initial states of interest, the associated total cost is near-optimal, \ie,
$
    |J(\hat{\bm{u}}(\,\cdot\,;0,\bm{x}_0);0,\bm{x}_0) - J(\bm{u}^*(\,\cdot\,;0,\bm{x}_0);0,\bm{x}_0)| \text{~is small.}
$

\subsection{Supervised-learning-based Approach for Closed-loop Optimal Control Problem}

Here we briefly explain the idea of the supervised-learning-based approach for the closed-loop optimal control problem.
The first step is to generate training data by solving the open-loop optimal control problems with zero initial time and initial states randomly sampled in $X$. Then, the training data is collected by  evenly choosing points in every optimal path:
\begin{equation*}
    \mathcal{D}=\{(t^{i,j},\bm{x}^{i,j}),\bm{u}^{i,j}\}_{1\le i \le M, 1\le j \le N},
\end{equation*}
where $M$ and $N$ are the number of sampled training trajectories and the number of points chosen in each path, respectively.
Finally, a function approximator (mostly neural network, as considered in this work) with parameters $\theta$ is trained by solving the following regression problem:
\begin{equation}\label{eq:NN_training}
    \min_{\theta} \frac{1}{MN}\sum_{i=1}^{M}\sum_{j=1}^{N}\|\bm{u}^{i,j} - \hat{\bm{u}}(t^{i,j},\bm{x}^{i,j};\theta)\|^2,
\end{equation}
and gives the Neural Network (NN) controller $\hat{\bm{u}}$. 

\revision{We emphasize that the basic supervised learning approach relies on two key assumptions. First, the open-loop optimal solutions must be unique, or at least the training data should contain only unimodal solutions. If multimodal open-loop solutions (e.g., a swing-up task that can proceed either left or right) are present in the training data, the function approximator will struggle to learn the closed-loop optimal control. Second, the function approximator must be capable of representing the closed-loop optimal control. According to the universal approximation theorem, neural networks can generally represent smooth control. However, standard neural networks may fail to represent discontinuous control, such as bang-bang control, which is beyond the scope of this work.}

\section{Proposed Method: Progressive Optimal Path Sampling}\label{Sec:POPS}
\begin{figure}[t]
     \centering
        \includegraphics[width=0.94\textwidth]{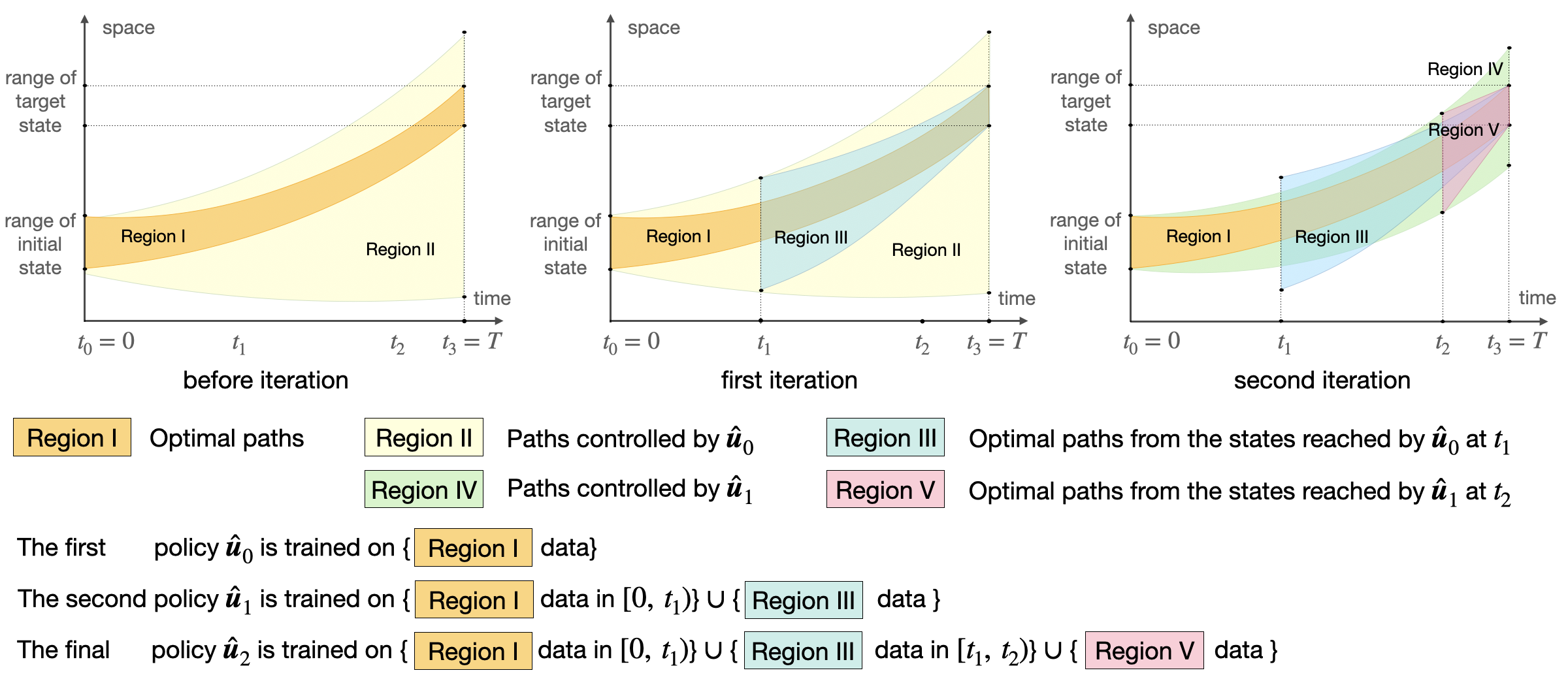}
        \caption{An illustration of POPS (Algorithm~\ref{alg:main}) when there are two intermediate temporal grid points $t_1$ and $t_2$.}
        \label{fig:ivp}
\end{figure}

Although the vanilla supervised-learning-based approach can  achieve a good performance in certain problems \citep{nakamura2021adaptive}, it is observed that its performance on complex problems is not satisfactory (see \citet{nakamura2021neural,zang2022machine} and examples below). One of the crucial reasons that the vanilla method fails is the distribution mismatch phenomenon. To better illustrate this phenomenon, let $\mu_0$ be the distribution of the initial state of interest and $\bm{u}:[0,T]\times \bR^n \rightarrow \mathcal{U}$ be a closed-loop control function. We use $\mu_{\bm{u}}(t)$ to denote the distribution of $\bm{x}(t)$ generated by $\bm{u}$:
$
    \dot{\bm{x}}(t) = \bm{f}(t,\bm{x}(t),\bm{u}(t,\bm{x}(t))), \bm{x}_0 \sim \mu_0.
$
Note that in the training process \eqref{eq:NN_training}, the distribution of the state at time $t$ is $\mu_{\bm{u}^*}(t)$, the state distribution generated by the closed-loop optimal control.
On the other hand, when we apply the learned NN controller in the dynamics, the distribution of the input state of $\hat{\bm{u}}$ at time $t$ is $\mu_{\hat{\bm{u}}}(t)$.
The error between state $\bm{x}$ driven by
$\bm{u}^*$ and $\hat{\bm{u}}$ accumulates and makes the discrepancy between $\mu_{\bm{u}^*}(t)$ and  $\mu_{\hat{\bm{u}}}(t)$ increases over time. Hence, the training data fails to represent the states encountered in the controlled process, and the error between $\bm{u}^*$ and $\hat{\bm{u}}$ dramatically increases when $t$ is large. See Figures \ref{fig:ivp} (left) and \ref{fig:sim2train} below for an illustration of this phenomenon.

\begin{algorithm}[ht] 
\caption{Progressive Optimal Path Sampling (POPS) for closed-loop optimal control design}
\label{alg:main}
\begin{algorithmic}[1]
\State \textbf{Input: }Initial distribution $\mu_0$, number of time points $K$, temporal grid points $0=t_0<t_1<\cdots<t_K=T$, time step $\delta$, number of initial points $N$. 
\State \textbf{Initialize: }  $S_{-1} = \emptyset$, $ \hat{\bm u}_{-1}(t,\bm x) = 0$.
\State Independently sample $N$ initial points from $\mu_0$ to get an initial point set $X_0$.

\For{$i=0,1,\cdots, K-1$}
\State For any $\bm x_0 \in  X_0$, compute IVP $(\bm x_0, 0, t_i, \hat{\bm  u}_{i-1})$
according to \eqref{eq:IVP}.  \Comment{Exploration}
\State Set $X_i = \{\bm x_{\hat{\bm u}_{i-1}}(t_i;0,x_0): x_0 \in X_0\}$.
\State For any $\bm x_i \in X_i$, call the open-loop optimal control solver to obtain $\bm x^*(t;t_i,\bm x_i) $ and $\bm u^*(t;t_i,\bm x_i)$ for $t \in [t_i,T]$.  \Comment{Labeling}
\State Set $\hat{S}_i = \{(t,\bm x^*(t;t_i,\bm x_i),\bm u^*(t;t_i,\bm x_i)): \bm x_i \in X_i, t\in [t_i,T], (t-t_0)/\delta \in \mathbb{N}\}.$
\State Set $S_i=\hat{S}_i\bigcup \{(t,\bm{x},\bm{u}): t<t_i,(t,\bm{x},\bm{u})\in S_{i-1} \}$.
 \State Train $\hat{\bm u}_{i}$ with dataset $S_i$. \Comment{Training}
\EndFor
\State \textbf{Output: } $\hat{\bm u}_{K-1}$.
\end{algorithmic}
\end{algorithm}

To overcome this problem, we propose the following progressive optimal path sampling (POPS) method. The key idea is to improve the quality of the NN controller iteratively by enlarging the training dataset with the states seen by the NN controller at previous times. Given predesigned (not necessarily even-spaced) temporal grid points $0 = t_0 < t_1 < \dots < t_K = T$, we first generate a training dataset $S_{0}$ by solving open-loop optimal control problems on the time interval $[0,T]$ starting from points in $X_0$, a set of initial points sampled from $\mu_0$, and train the initial model $\hat {\bm u}_{0}$. Under the control of $\hat {\bm u}_0$, the generated trajectory deviates more and more from the optimal trajectory. So we stop at time $t_1$, \ie, compute the IVPs~\eqref{eq:IVP} using $\hat {\bm u}_0$ as the closed-loop control and points in $X_0$ as the initial points on the time interval $[0,t_1]$, and then on the interval $[t_1,T]$ solve new optimal paths that start at the end-points of the previous IVPs. 
The new training dataset $S_1$ is then composed of new data (between $t_1$ and T) and the data before time $t_1$ in the dataset $S_0$, and we train a new model $\hat{ \bm{u}}_1$ using $S_1$. We repeat this process to the predesigned temporal grid points $t_2, t_3, \cdots$ until end up with $T$. 
In other words, in each iteration, the adaptively sampled data replaces the corresponding data (defined on the same time interval) in the training dataset (the size of the training data remains the same). The whole process can be formulated as Algorithm~\ref{alg:main}, and we refer to Figure \ref{fig:ivp} for an illustration of the algorithm's mechanism.
We refer to this method as \textit{progressive optimal path sampling} because it incrementally generates new open-loop optimal paths as training data using the most recent NN controller. While this approach shares similarities with a few existing methods, as discussed in the next section, our strategy for progressively adding data is more appropriate for optimal control compared to the aggressive selection used in methods such as \dagger.

It is worthwhile mentioning that POPS is versatile enough to combine other improvements for closed-loop optimal control problems, such as efficient open-loop control problem solvers \citep{kang2021algorithms,zang2022machine} or specialized neural network structures \citep{qrnet1,nakamura2021neural,nakamura2022neural}. One design choice regarding the network structure in POPS is whether to share the same network among different time intervals. We choose to use the same network for all the time intervals in the following numerical examples, but the opposite choice is also feasible.

\section{Related Work}
\label{sec:compare}
\paragraph{Adaptive sampling for supervised-learning-based approach.}
\revision{
Learning the optimal policy can be framed as an imitation learning problem by treating the optimal policy as the expert's control. This leads to a similar distribution mismatch issue in both settings. However, there is a key difference regarding the mechanism of data generation between the two settings: in imitation learning, it is often assumed that one can easily access the expert's behavior at every time-state pair, \eg, through solving a sub-problem by MPC~\citep{PLATO2016,APEX2021}, while in the optimal control problem, it is much more computationally expensive to access since one must solve an open-loop optimal control problem.
This difference affects algorithm design fundamentally. 
The methods for imitation learning often assume low computational costs of obtaining the expert demonstration at each time-state pair.
Take the forward training algorithm~\citep{ross2010efficient} as an example.
To apply it to the closed-loop optimal control problem, we first need to consider a discrete-time version of the problem with a sufficiently fine time grid: $0=t_0 < t_1<\dots<t_{K'} = T$. At each time step $t_i$, we learn a policy function $\bar{\bm{u}}^i:\bR^n \rightarrow \mathcal{U}$ where the state $x$ in the training data are generated by sequentially applying $\bar{\bm{u}}^0,\dots,\bar{\bm{u}}^{i-1}$ and the labels are generated by solving the open-loop optimal solutions with $(t_i, x)$ as the initial time-state pair. Hence, the open-loop control solver is called with numbers proportionally to the discretized time steps $K'$, and only the first value on each optimal control path is used for learning. In contrast, in Algorithm \ref{alg:main}, we can use much more values over the optimal control paths in learning, which allows a much coarse temporal grid for adaptive sampling, and significantly reduces the total cost of solving open-loop optimal control problems.}

\revision{
Another popular imitation learning method is \dagger (Dataset Aggregation) \citep{ross2011reduction}, which can be adapted for sampling in closed-loop optimal control. Similar to the forward training algorithm, \dagger assumes access to expert behavior at each time-state pair. Here, we retain the core idea of \dagger but modify it for closed-loop optimal control.
In \dagger, in order to improve the current closed-loop controller $\hat{\bm{u}}$, one simulates the system using $\hat{\bm{u}}$ over $[0, T]$ starting from various initial states and collect the states on a time grid $0=t_0  < t_1 < \dots < t_{K-1} < T$. The open-loop control problems are then solved with \textit{all} the collected time-state pairs as the initial time-state pair, and \textit{all} the corresponding optimal solutions are used to construct a dataset for learning a new controller. The process can be repeated until a good controller is obtained.
While \dagger aims to address distribution mismatch, its state selection differs from POPS. Take the data collection using the controller $\hat{\bm{u}}_1$ in the first iterative step for example. POPS focuses on the states at the time grid $t_1$ while \dagger collects states at all the time grids. If $\hat{\bm{u}}_1$ is still far from optimal, the data collected at later time grids may be irrelevant or even misleading due to error accumulation in states. As shown in subsequent sections, this sensitivity can reduce \dagger's effectiveness compared to POPS.}

There are other adaptive sampling methods specifically developed for closed-loop optimal control rather than imitation learning.  \citet{VariationalPolicySearchTO2013,levine2014learning} propose iterative methods that generate sub-optimal open-loop solutions around the current policy but are tailored to stochastic policies, limiting their applicability to deterministic control problems considered in this paper. \citet{nakamura2021adaptive} propose an adaptive sampling method that chooses the initial points with large gradients of the value function as the value function tends to be steep and hard to learn around these points. 
\citet{landry21seagul} propose to sample the initial points on which the errors between predicted values from the NN and optimal values are large. 
These two adaptive sampling methods both focus on finding points that are not learned well but ignore the influence of the accumulation of the distribution mismatch over time brought by controlled dynamics.

\paragraph{Open-loop data for direct policy search approach.}
\revision{
Open-loop data can also alleviate optimization challenges in direct policy search. For example, guided policy search~\citep{GPS2013} incorporates trajectory samples from differentiable dynamic programming into the optimization objective. Similarly, \citet{combiningTO2014} reformulate the optimization as a constrained problem to prevent the learned policy from deviating too far from open-loop optimal trajectories.}

\section{Theoretical Analysis on an LQR Example}\label{sec:toy_example}
In this section, we analyze the superiority of POPS by considering the following uni-dimensional linear quadratic regulator (LQR) problem:
\begin{align*}
    &\min_{x(t),u(t)} \frac{1}{T}\int_{t_0}^{T}|u(t)|^2\rmd t + |x(T)|^2 \\
    \text{s.t. }& 
    \dot{x}(t) = u(t), \, t \in [t_0,T], \quad
    x(t_0) = x_0,
\end{align*}
where $T$ is a positive integer, $t_0 \in [0,T]$ and $x_0 \in \bR$. Classical theory on linear quadratic control (see, \eg \citet{sontag2013mathematical}) gives the following explicit linear form of the optimal controls:
\begin{numcases}{\hspace{-4em}}
    u^*(t;t_0,x_0) = \textstyle{-\frac{T}{T(T-t_0) + 1}x_0}, 
    \tag*{(open-loop optimal control)} \\ 
    u^*(t,x)= \textstyle{-\frac{T}{T(T-t)+1}x}.  \tag*{(closed-loop optimal control)}
\end{numcases}
We consider the following two models to approximate the closed-loop optimal control function with parameter $\theta$:
\begin{align}
    &\text{Model 1:~~~~}u_\theta(t,x) = \textstyle{-\frac{T}{T(T-t)+1}x} + b(t), \text{ where } \theta =\{\theta_t\}_{0\le t\le T} = \{b(t)\}_{0 \le t \le T}. \label{lqr_model1} \\
    &\text{Model 2:~~~~}u_\theta(t,x) = a(t) x + b(t), \text{ where } \theta = \{\theta_t\}_{0\le t\le T}= \{(a(t),b(t))\}_{0\le t \le T}. \label{lqr_model2}
\end{align}
Since learning a linear model has no error with exact data, to mimic the errors encountered when learning neural networks, throughout this section, we assume the data has certain noise. 
To be precise, for any $t_0 \in [0,T]$ and $x_0 \in \bR$, the open-loop optimal control solver gives the following approximated optimal path:
\begin{align*}
\begin{cases}
    \hat{u}(t;t_0,x_0) = -\frac{T}{T(T-t_0) + 1}x_0 + \epsilon Z, \vspace{0.5em} \\ 
    \hat{x}(t;t_0,x_0) = x_0 +\int_{t_0}^{t} \hat{u}(t;t_0,x_0)\rmd t = \frac{T(T-t)+1}{T(T-t_0) +1}x_0 + (t-t_0)\epsilon Z,
\end{cases}
\end{align*}
where $\epsilon > 0$ is a small positive number to indicate the scale of the error and $Z$ is a normal random variable whose mean is $m$ and variance is $\sigma^2$. In other words, the obtained open-loop control is still constant in each path, just like the optimal open-loop control, but perturbed by a random constant. The random variables in different approximated optimal paths starting from different $t_0$ or $x_0$ are assumed to be independent.

\paragraph{Comparison results.} 
\revision{
We provide a theoretical comparison under Model 1~\eqref{lqr_model1} and a numerical evaluation under Model 2~\eqref{lqr_model2} for the vanilla supervised-learning-based method, \dagger, and POPS. Detailed descriptions of the setups of each method are given in Appendix~\ref{appendix_lqr}. Notably, we maintain the same total number of open-loop optimal paths, $NT$, across all methods, where $N$ is an integer.}

We denote four closed-loop controllers as $u_o$, $u_v$, $u_d$, and $u_p$, corresponding to the optimal controller, the controller learned by the vanilla method, the controller learned by \dagger, and the controller learned by POPS, respectively. Additionally, we define $\hat{x}^v(t)$, $\hat{x}^d(t)$, and $\hat{x}^p(t)$ as random variables whose distributions are the average distributions of state variables in the training data at time $t$ for the vanilla method, \dagger, and POPS (in the final iteration), respectively. 
We define two key metrics for evaluation: (1) \textit{distribution difference}, which measures the discrepancy between the variance of the training trajectory data used for the final controller and the variance of trajectory data generated by applying the final controller, and (2) \textit{performance difference}, which is the difference between the total cost incurred by the final learned controller and the optimal cost. Theorem~\ref{thm_1} establishes that, for both the vanilla method and \dagger, the distribution difference and performance difference increase quadratically with $T$, while in the case of POPS, these metrics remain independent of $T$.
Therefore, compared to the vanilla method and \dagger, POPS mitigates the distribution mismatch phenomenon and significantly improves the performance for large $T$. The detailed proof is provided in Appendix~\ref{appendix_lqr}.

\begin{thm}\label{thm_1}
Under Model 1~\eqref{lqr_model1}, 
define dynamical systems:
$\displaystyle{
				\dot{x}_s(t) = u_s(t) = u_s(t,x_s(t)),\, x_s(0) = \xinit}
$, ${0 \le t \le T,\, s \in \{o,v,d,p\}}.$
\begin{enumerate}
    \item \textbf{Distribution difference:}
    \revision{Assume $\xinit$ is a random variable following a standard normal distribution, which is independent of the initial points and noises in the training process. Then, $\displaystyle{\mathbb{E} \hat{x}^v(t) = \mathbb{E}x_v(t), \mathbb{E}\hat{x}^p(t) = \mathbb{E}x_p(t)}$ and
    \begin{align*}
        &|\Var(\hat{x}^v(t)) - \Var(x_v(t))| = \epsilon^2\sigma^2(1-\frac{1}{NT}) t^2, \\
        &\Var(\hat{x}^d(t)) - \Var(x_d(t)) \ge \epsilon^2\sigma^2(\frac{t^2}{3}-\frac{2t}{N}), \\
        &|\Var(\hat{x}^p(t)) - \Var(x_p(t))| \le \epsilon^2\sigma^2.
    \end{align*}
    }
    \item \textbf{Performance difference:} Assume $\xinit$ is fixed and define the total cost
    $
        J_s =\frac{1}{T} \int_{0}^{T}|u_s(t)|^2\rmd t + |x_s(T)|^2, s \in \{o,v,d,p\}.
    $
    Then,
\begin{align*}
       &\mathbb{E} J_v - J_o =(T^2+1)(m^2+\frac{\sigma^2}{NT})\epsilon^2 , \\
       &\mathbb{E}J_d - J_o \ge (\frac{T^2m^2}{4}+\frac{T\sigma^2}{3N})\epsilon^2, \\
       &\mathbb{E} J_p - J_o \le 3(m^2+\frac{\sigma^2}{N})\epsilon^2. 
\end{align*}

\end{enumerate}
\end{thm}

\begin{figure}[t]
        \centering
            \includegraphics[width=0.45\textwidth]{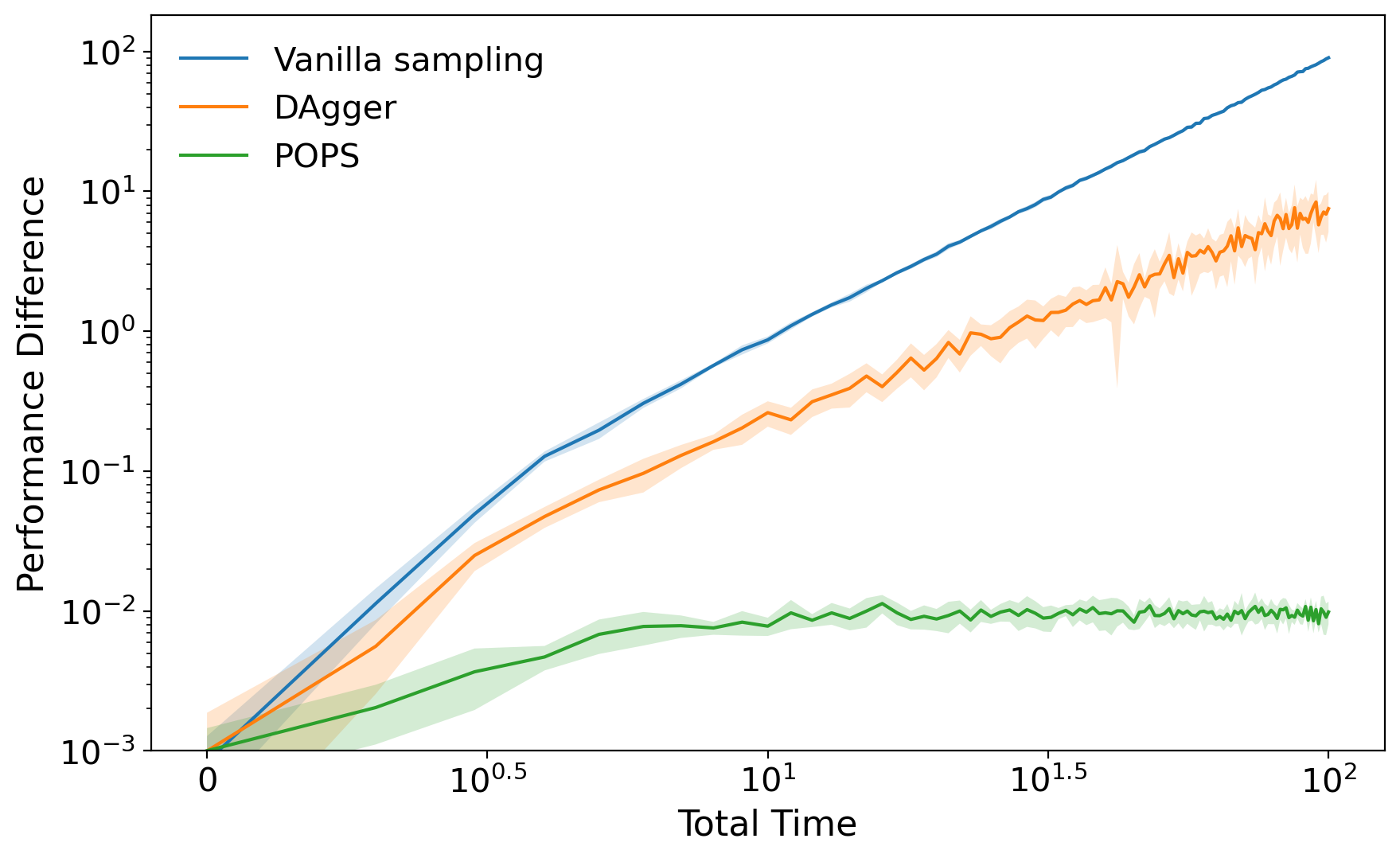}  \quad\quad
            \includegraphics[width=0.45\textwidth]{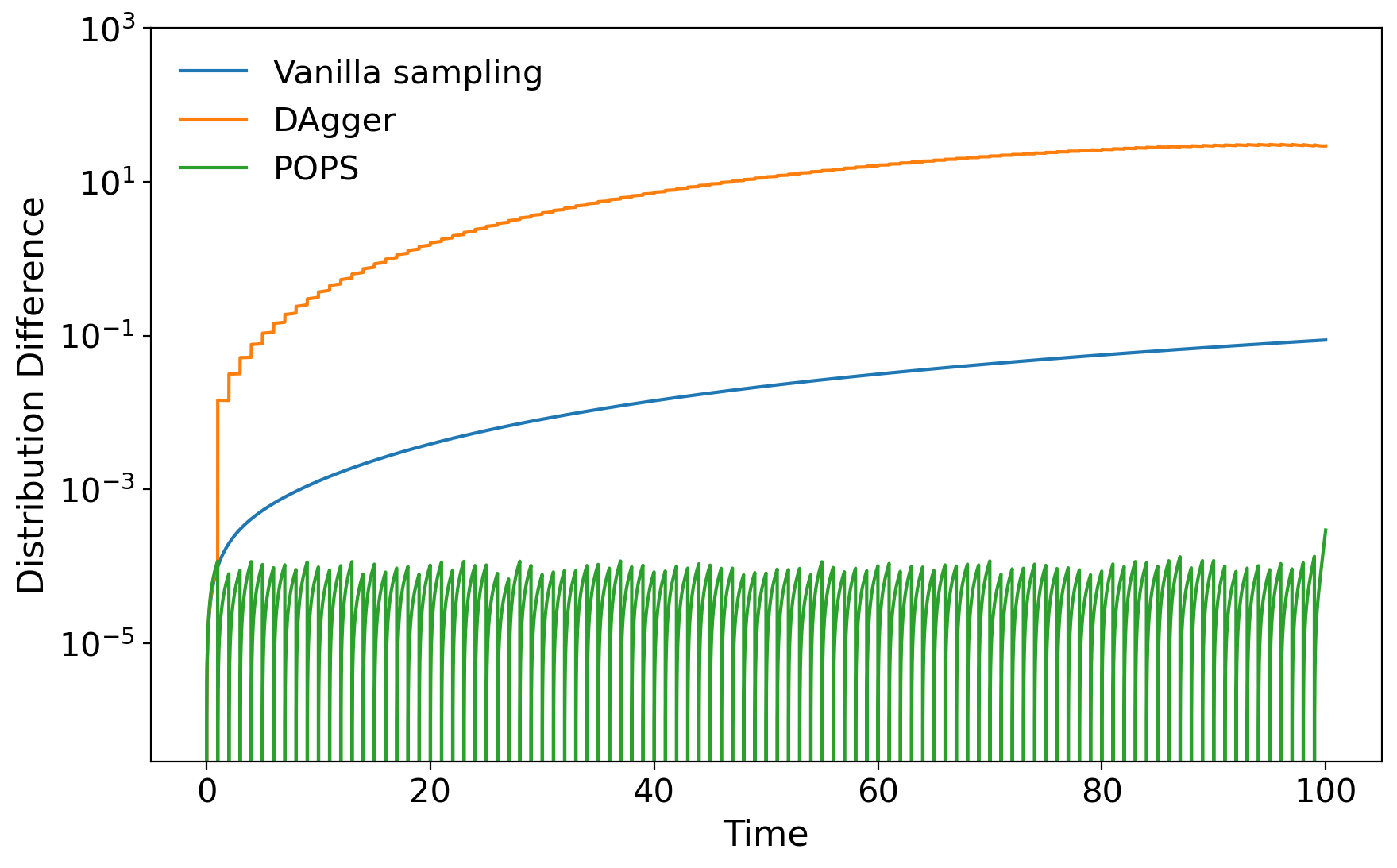}\\
            \includegraphics[width=0.45\textwidth]{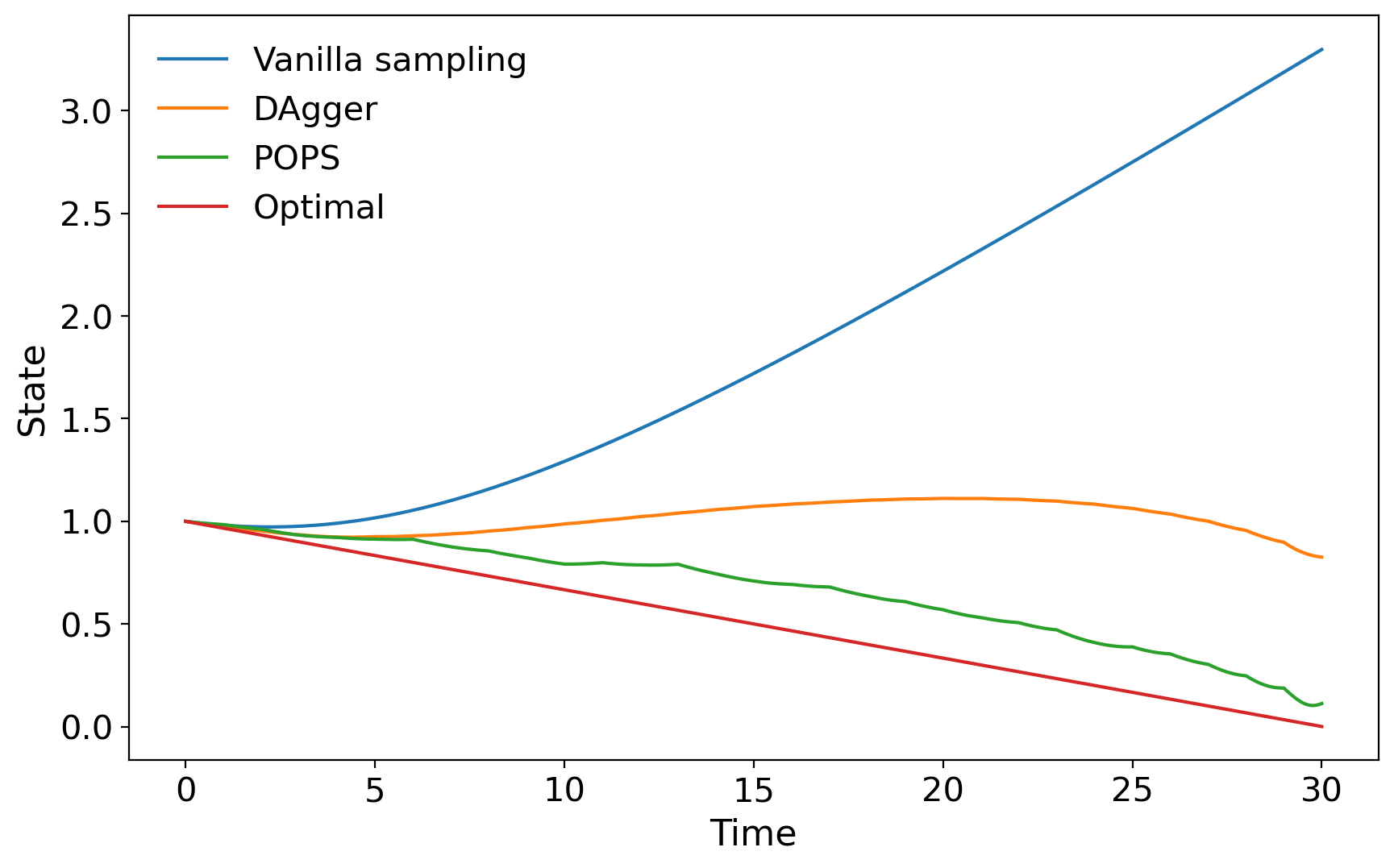} \quad\quad
            \includegraphics[width=0.45\textwidth]{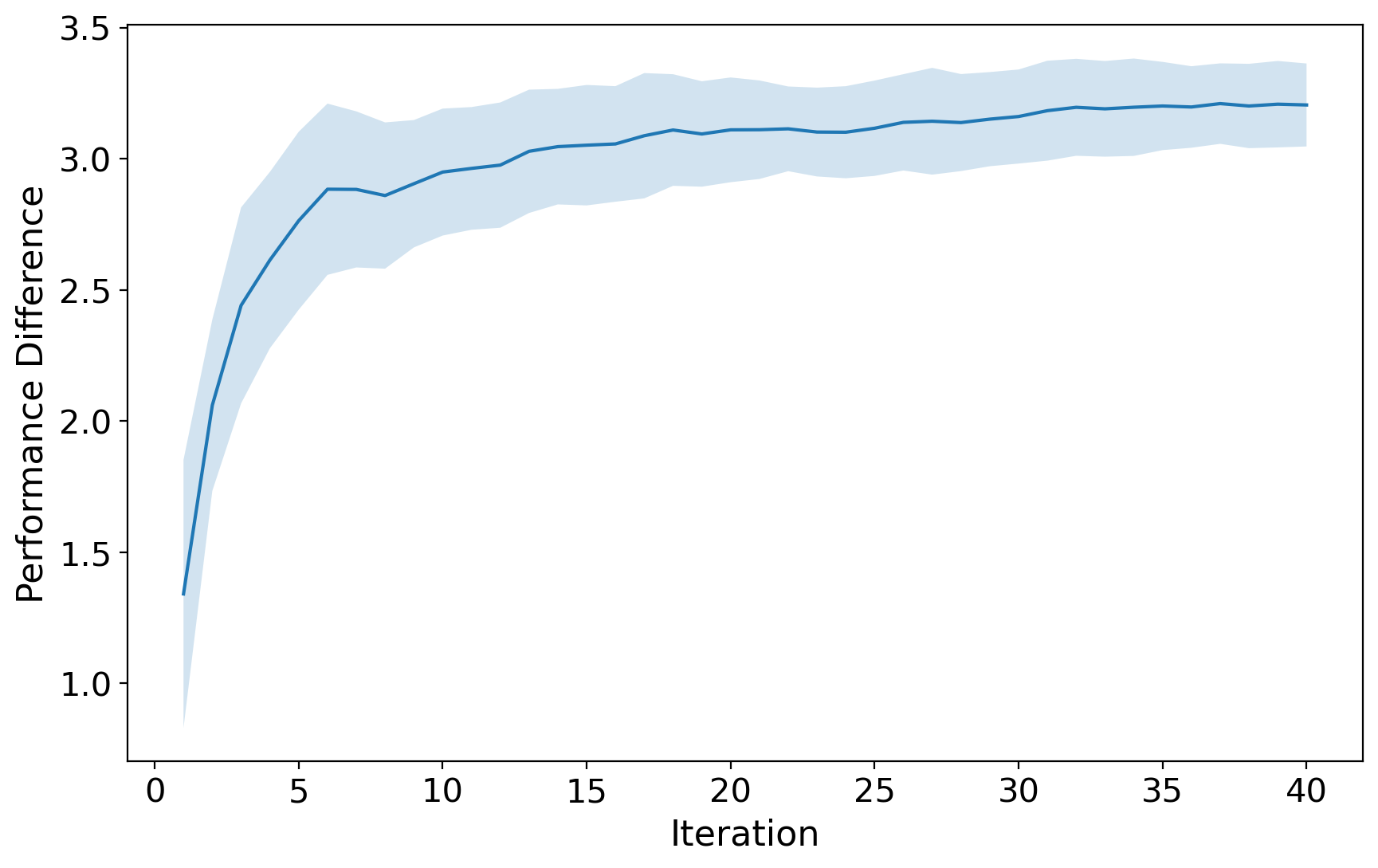}
         \caption{\revision{Numerical results on learning Model 2~\eqref{lqr_model2} for the LQR model. \textbf{Top left:}  performance differences (in the logarithm scale) of the vanilla sampling method, \dagger, and POPS for different total time span $T$ (in the logarithm scale). \textbf{Top right:} differences of the second-order moments (in the logarithm scale) between the distributions of the training data and the data reached by the controllers at intermediate times. \textbf{Bottom left:} the optimal path and the paths generated by the vanilla sampling method, \dagger, and POPS. \textbf{Bottom right:} the performance difference of \dagger along multiple iterations.}}
        \label{fig:lqr}
    \end{figure}

\revision{
Next, we present the numerical results when we use Model 2~\eqref{lqr_model2} to fit the closed-loop optimal control, with the setting $\epsilon = 0.1$, $m = 0.1$ and $\sigma^2 = 1$.
Figure \ref{fig:lqr} (top left) compares the performance of the vanilla method, \dagger, and POPS on different total times $T$.  The performance difference is an empirical estimation of $\mathbb{E} [J_v - J_o]$, $\mathbb{E} [J_d - J_o]$ and $\mathbb{E}[J_p - J_o]$ when $\xinit$ follows a standard normal distribution. In this experiment, for each method, we set $N = 100$ and learn 10 different controllers with different realizations of the training data and calculate the average of the performance difference on 1000 randomly sampled initial points (from a standard normal distribution) and 10 learned controllers. Figure \ref{fig:lqr} (top right) shows how the time $t$ influences the distribution differences when $T = 100$ and $N = 100$.
Figure \ref{fig:lqr} (bottom left) compares the optimal path $x_o(t)$ with $x_v(t)$, $x_d(t)$ and $x_p(t)$, the trajectory generated by the controllers learned by the vanilla method, \dagger, and POPS starting from $\xinit = 1$, when $T = 30$ and $N = 100$. All three experiments consistently demonstrate the superior performance of POPS as similarly established in Theorem~\ref{thm_1}. Finally, Figure \ref{fig:lqr} (bottom right) shows the performance of \dagger with multiple iterations under $T = 30$. The results indicate that increasing the number of iterations does not improve the controller’s performance in \dagger.
}

\revisionjmlb{We conclude this section by noting that, although the analysis here is specific to linear-quadratic settings, it illustrates how POPS effectively controls the growth of distribution mismatch with respect to the time horizon. Unlike single-step imitation approaches such as DAgger, POPS leverages entire open-loop trajectories, which amortizes the cost of data generation from optimal control solvers and helps mitigate long-horizon distribution shift. This is supported both by the theoretical results in this section and by the nonlinear examples presented in the next two sections.
}

\section{The Optimal Landing Problem of Quadrotor}\label{sec:prob_quadrotor}
\revisionjmla{In this section, we evaluate the performance of POPS on a quadrotor optimal-landing task: bringing the vehicle to a prescribed touchdown point while minimizing energy consumption.} We consider the full quadrotor dynamic model with  12-dimensional state variable and 4-dimensional control variable~\citep{bouabdallah2004design,madani2006control,mahony2012multirotor}. The state variable is represented as $\bm{x}=(\bm{p}\transpose,\bm{v}_b\transpose,\bm{\eta}\transpose, \bm{w}_b\transpose)\transpose \in \mathbb{R}^{12}$, where $\bm{p}=(x,y,z)\in\mathbb{R}^3$ denotes the quadrotor's position in Earth-fixed coordinates, $\bm{v}_b\in\mathbb{R}^3$ is the velocity in body-fixed coordinates, $\bm{\eta}=(\phi, \theta, \psi) \in \mathbb{R}^3$ represents the attitude (roll, pitch, yaw) in Earth-fixed coordinates, and $\bm{w}_b \in \mathbb{R}^3$ is the angular velocity in body-fixed coordinates. The control variable is defined as $\bm{u} = (s, \tau_x, \tau_y, \tau_z)\transpose \in \mathbb{R}^4$, where $s$ is the total thrust, and $\tau_x, \tau_y, \tau_z$ are the body torques generated by the four rotors. The details of the quadrotor dynamics are provided in Appendix~\ref{appendix:dyn}. 

Our goal is to compute optimal landing trajectories from an initial state $\bm{x}_0$ to a target state $\bm{x}_T = 0$ with minimum control efforts over a fixed time span $T=16$. The initial state distribution is uniform over the set 
$$ X = \{x, y \in [-40, 40], z \in [20, 40], v_x, v_y, v_z \in [-1, 1], \theta, \phi \in [-\pi/4, \pi/4], \psi \in [-\pi, \pi]; \bm{w}_b = \bm{0}\}.
\label{eq:quad_init_set}$$ 
The running cost and terminal cost are defined as follows:
\[
L(\bm{x}, \bm{u}) = (\bm{u} - \bm{u}_d)\transpose Q_u (\bm{u} - \bm{u}_d), \quad
M(\bm{x}) = \bm{p}\transpose Q_{pf} \bm{p} + \bm{v}\transpose Q_{vf} \bm{v} + \bm{\eta}\transpose Q_{\eta f} \bm{\eta} + \bm{w}\transpose Q_{wf} \bm{w} = \bm{x}\transpose Q_f \bm{x},
\]
where $\bm{u}_d = (mg, 0, 0, 0)$ is the reference control to counteract gravity, and $Q_u = \text{diag}(1, 1, 1, 1)$ is the weight matrix that penalizes deviations from the reference control. The weight matrices for the terminal cost are $Q_{pf} = 5I_3$, $Q_{vf} = 10I_3$, $Q_{\eta f} = 25I_3$, and $Q_{wf} = 50I_3$. Larger weights are used in the terminal cost to impose stricter penalties on deviations from the landing target. The open-loop optimal solutions are derived by solving the associated two-point boundary value problem using the space-marching technique~\citep{zang2022machine}; additional details are provided in Appendix~\ref{appendix:PMP}. When applying POPS, the number of initial points for these open-loop optimal solutions is fixed at $N = 500$ for each iteration. To construct the training dataset, we sample time-state-action tuples at intervals of $\delta = 0.2$ time steps along the optimal trajectories, resulting in a consistent dataset size of $81 \times 500$ at every iteration.

The neural network models used in all quadrotor experiments share the same structure: a 13-dimensional input (12 for the state variables and 1 for time) and a 4-dimensional output. Each network is fully connected, with two hidden layers comprising 128 neurons per layer. We use the \texttt{tanh} activation function in both layers. The inputs are scaled to the range $[-1, 1]$, where the upper and lower bounds correspond to the maximum and minimum values of the training dataset. Since the activation function is \texttt{tanh}, we apply Xavier initialization~\citep{Xavier_ini} prior to training. The neural networks are trained using the Adam optimizer~\citep{adam} with a learning rate of 0.001, a batch size of 1000, and 1000 epochs. At each iteration of POPS, a new neural network is trained from scratch. The quadrotor experiments were conducted on a MacBook Pro equipped with an Apple M1 Pro chip.

\begin{figure}[t]
\centering
\includegraphics[width=0.9\textwidth]{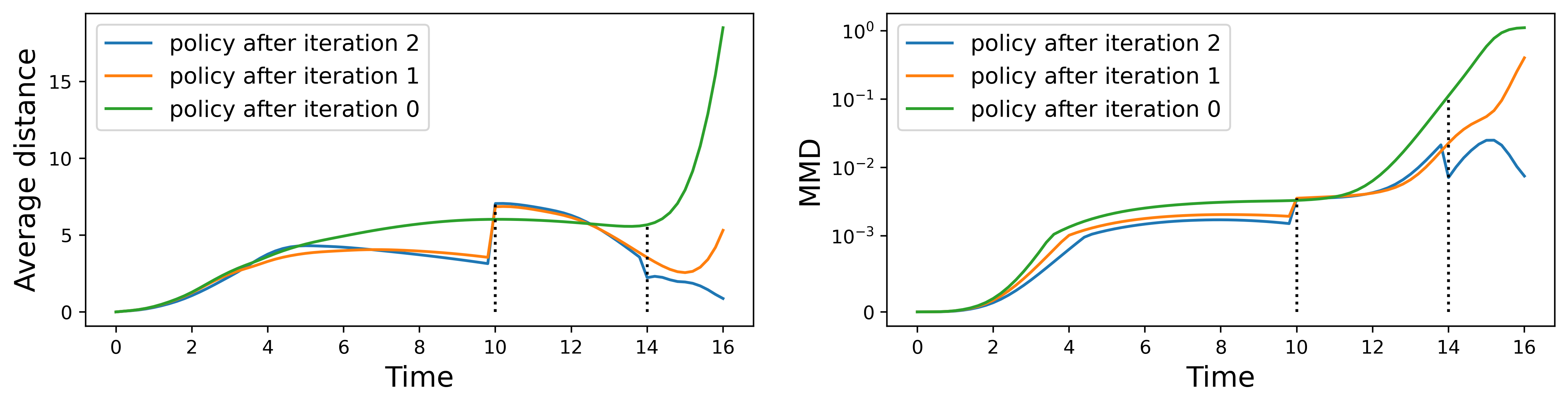}
\caption{Left: the average pointwise distance between the training data and the data reached by controllers at different times. Right: the maximum mean discrepancy (in the logarithm scale) between the training data and the data reached by controllers at every time using the Gaussian kernel.}
\label{fig:sim2train}
\end{figure}

In the first experiment, we apply POPS and select the temporal grid points $t = 0, 10, 14, 16$. After training, we use the learned NN controllers to solve the initial value problem for the 500 training initial points, comparing the trajectories driven by the NN controllers to their corresponding training data. In Figure~\ref{fig:sim2train}, the left sub-figure displays the average distance between the states reached by the NN controller and the states from the training data at various time steps. The right sub-figure shows the maximum mean discrepancy (MMD)~\citep{borgwardt2006integrating} between these two datasets, using a Gaussian kernel $k(x, y) = \exp(-\frac{|x - y|^2}{2})$.
Both figures show noticeable jumps at $t = 10$ and $14$ because the NN-controlled trajectory is continuous across time, whereas the training data is discontinuous at the points where POPS is applied. Notably, without adaptive sampling (represented by the curve labeled ``policy after iteration 0" in Figure~\ref{fig:sim2train}), there is a substantial discrepancy between the NN-controlled states and the training data while this discrepancy progressively decreases with each iteration of POPS.

\begin{figure}[t]
\centering
\includegraphics[width=0.95\textwidth]{fig/review_traj.png}
 \caption{The optimal path and path controlled by learned controllers. We show the 3-dimensional position $\bm p=(x,y,z)$ and 3-dimensional attitude $\bm{\eta}=(\phi,\theta,\psi)$ in terms of Euler angles in Earth-fixed coordinates. The cost of $\hat{u}_0, \hat{u}_1,\hat{u}_2$ controlled paths is 3296.2, 119.9, 6.7, respectively, and the optimal cost is 6.3.}
\label{fig:quad_traj}
\end{figure}

Next, we report the performance of the learned NN controllers. Example trajectories driven by NN controllers are illustrated in Figure \ref{fig:quad_traj}. 
As the POPS iterations progress, the trajectories increasingly align with the optimal path. In the final iteration, the trajectory controlled by $\hat{\bm u}_3$ closely matches the entire optimal trajectory. Quantitatively, the costs of the three controlled trajectories are 3296.2, 119.9, and 6.7, respectively, compared to the optimal cost of 6.3.

\begin{figure}[t]
\begin{floatrow}
\ffigbox[\FBwidth]{%
  \includegraphics[width=0.34\textwidth,height =0.25\textwidth]{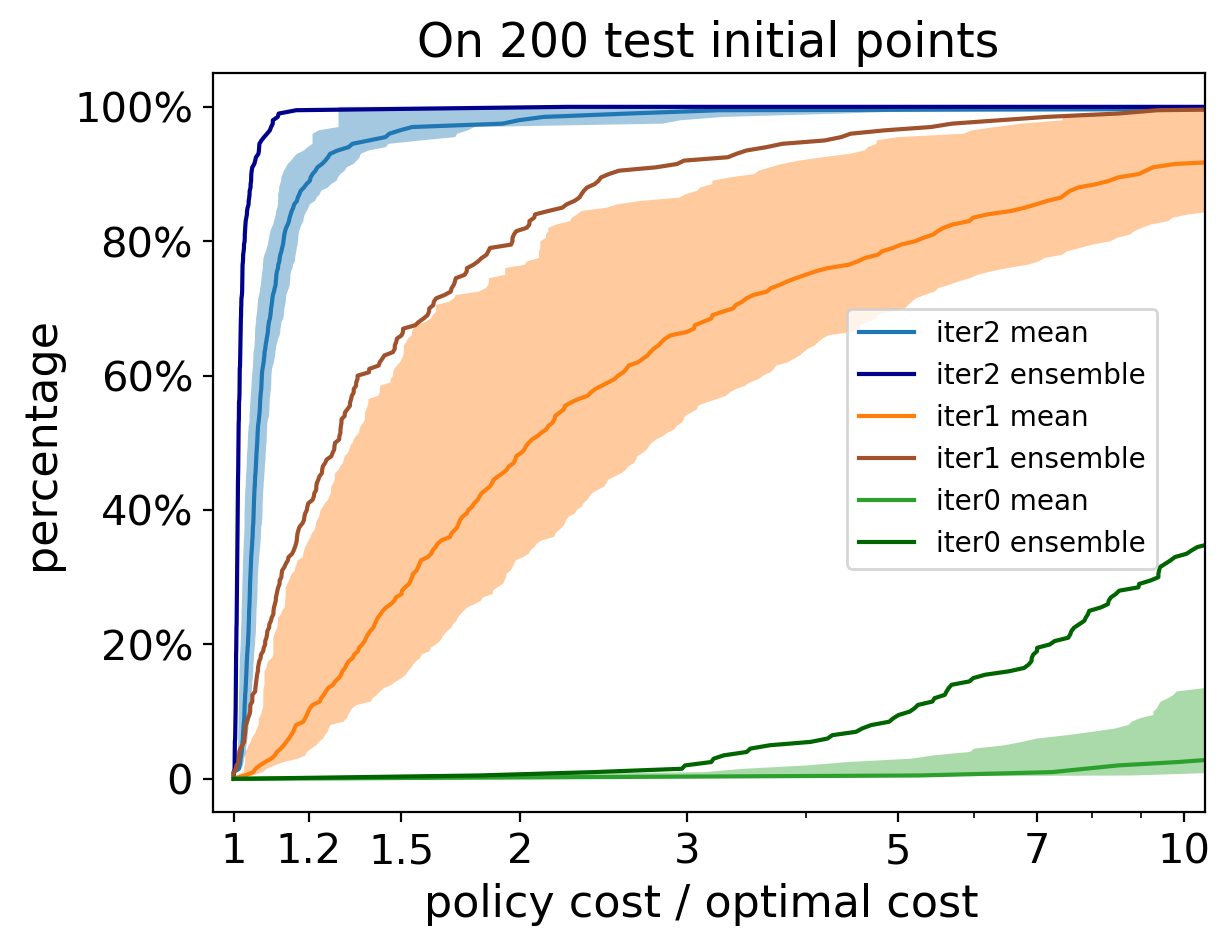}%
}{%
  \caption{Cumulative distribution function of the cost ratio  %
   of POPS \revision{among 200 test trajectories}.
   The shaded area represents the mean cost ratio $\pm$ the standard deviation.
   }
   \label{fig:quadrotor_200fig}
}
\capbtabbox{%
\footnotesize
  \begin{tabular}{cccc}
    \toprule
\textbf{Policy}  & \textbf{Mean}   & \textbf{90\%}  & \textbf{Median} \\ \midrule
\textbf{$\hat{\bm u}_0$} &   16.43 ($\pm$ 6.16) & 17.60 ($\pm$ 4.80)  & 16.54 ($\pm$ 6.92)  \\
\textbf{$\hat{\bm u}_1$}  &  3.69 ($\pm$ 1.62) & 9.11 ($\pm$ 5.54) & 2.06 ($\pm$ 0.72)  \\
\textbf{$\hat{\bm u}_2$} &   1.17 ($\pm$ 0.09) & 1.22 ($\pm$ 0.09) &  1.06 ($\pm$ 0.02)  \\
\textbf{$\hat{\bm u}_0$ ensemble} & $20.37$ & $46.85$ & $15.05$ \\
\textbf{$\hat{\bm u}_1$ ensemble} & $1.78$ & $2.54$ &$1.29$ \\
\textbf{$\hat{\bm u}_2$ ensemble} & $1.03$ & $1.04$ & $1.01$ \\
\bottomrule
    \end{tabular}
}{%
  \caption{\revision{The columns ``Mean'' ``90\%'', and ``Median'' represent the cost ratio statistics across 200 test initial points.} The first 3 rows for $\hat{\bm u}_0$, $\hat{\bm u}_1$, $\hat{\bm u}_2$ denote the policy after the first, second, and third round of training, respectively and the entries are averages and standard deviations of 5 experiments. The bottom 3 rows report the results from the ensemble of five networks trained independently. }
    \label{fig:quadrotor_200table}
}
\end{floatrow}
\end{figure}

\revision{
To evaluate the controller's performance, we use the cost ratio, defined as the total cost \eqref{eq:loss} under the given controller divided by the optimal cost. By definition, the cost ratio is always greater than or equal to 1.  We perform experiments five times with different random seeds and present the results in Figure~\ref{fig:quadrotor_200fig} and Table~\ref{fig:quadrotor_200table}. The performance is evaluated on 200 test initial points which are also uniformly drawn from the initial set $X$ defined in \eqref{eq:quad_init_set}. Figure~\ref{fig:quadrotor_200fig} shows the cumulative distribution function (CDF) of the cost ratio under various NN controllers on test initial points. The CDF of the optimal controller’s cost ratio is a horizontal line at a ratio of 1 and a percentage of 100\%. Therefore, a curve closer to the top-left corner indicates better performance of the corresponding controller. It is important to note that the randomness in the first iteration arises from the initialization and batch data sampling during the training of the neural networks. In subsequent iterations, additional randomness is introduced through both the generation of distinct training data by different NN controllers and the stochastic nature of their training processes. The results indicate that the POPS algorithm improves performance with each iteration, reducing the cost from the initial model $\hat{\bm{u}}_0$, which has an average cost ratio of 16.43, to $\hat{\bm{u}}_2$, which has an average cost ratio of 1.17. Additionally, we evaluate ensembles that average the output of five independently trained networks (indexed by $i$) as a closed-loop controller: $\hat{\bm u}^{{\textbf{ensemble}}}(t, x) = \frac{1}{5} \sum_{i=1}^5 \hat{\bm u}_{i}(t, x)$, represented by the dark curves. This ensemble demonstrates superior performance over individual networks, achieving an average cost ratio of 1.03.
}

\begin{figure}[t]
\centering
\includegraphics[width=0.95\textwidth]{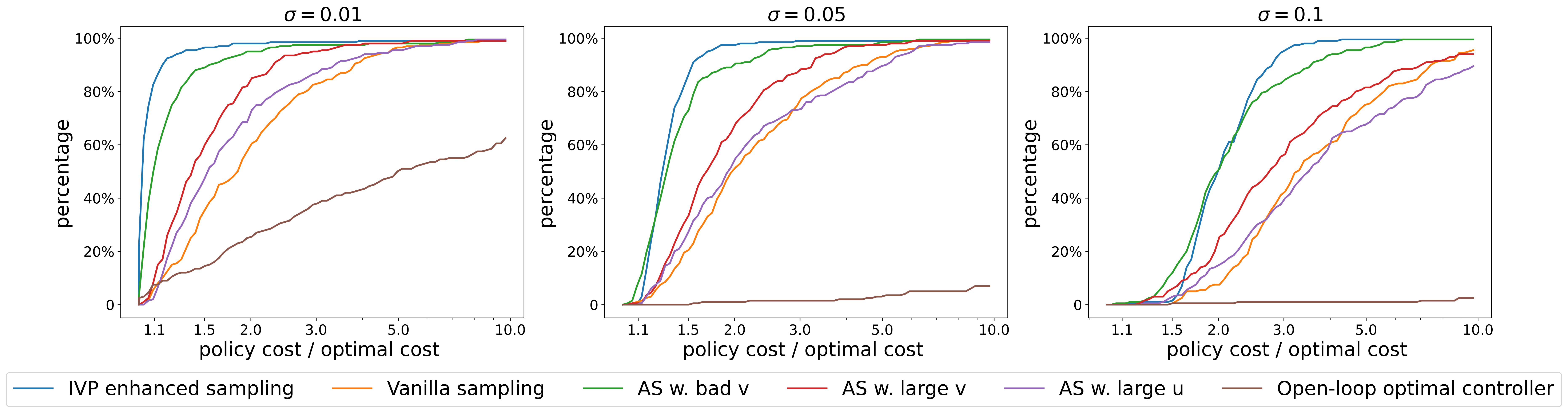}
 \caption{Cumulative distribution function of the cost ratio between NN controlled value under disturbance and the optimal value \revision{among 200 test trajectories}.}
\label{fig:dt}
\end{figure}

We further assess the performance of NN controllers under observation noise, simulating real-world sensor errors.  During simulation, we add a disturbance $\epsilon$ (including time) uniformly sampled from $[-\sigma,\sigma]^{13}$ to the network input.
We test with $\sigma=0.01,0.05,0.1$, and report the results in Figure \ref{fig:dt}. For comparison, we also evaluate the performance of the open-loop optimal controller with a disturbance $\hat \epsilon \in \mathbb{R}$ added to the input time.
Figure \ref{fig:dt} demonstrates that closed-loop controllers are more robust than open-loop controllers, and the controller trained with POPS continues to perform well under these noise levels.

We also test four different temporal grid point settings in Algorithm~\ref{alg:main} to assess their impact. The results in Table \ref{tab:quadrotor_time_grid} demonstrate the algorithm’s robustness, with consistent performance across different choices of temporal grid points.

\begin{table}[!htbp]
\centering
\begin{tabular}{c|ccccc}
\toprule
&\multicolumn{5}{c}{\textbf{Temporal grid points ($T=16$)}}\\
\textbf{$\#$ of iterations}&4&8&10&12&14\\
\midrule
2&&&&&  1.36 ($\pm$ 0.08) \\
3& &&  3.69 ($\pm$ 1.62) &&  1.17 ($\pm$ 0.09) \\
4& &  10.86 ($\pm$ 2.27) &&  1.72 ($\pm$ 0.49) &  1.17 ($\pm$ 0.07) \\
6&   14.39 ($\pm$ 1.85) &  12.14 ($\pm$ 5.64) &  6.30 ($\pm$ 2.67) &  1.73 ($\pm$ 0.46) &  1.22 ($\pm$ 0.04)\\
\bottomrule
\end{tabular}
\caption{Average cost ratio on 200 test points of POPS applied to different temporal grid points. The outcomes were consistent after the initial iteration (16.43 ($\pm$ 6.16)), and thus, we have excluded them from the table. This uniformity arises from utilizing the same five random seeds for repeating the experiments. The first line indicates that we performed 2 iterations, and the corresponding temporal grid points for adaptive sampling are $0<14<16$. After iteration 0, the average ratio of policy cost over optimal cost is 16.43 and after iteration 1, it decreases to 1.36. The second line shows the same experiment in Figure \ref{fig:quadrotor_200fig}.}
\label{tab:quadrotor_time_grid}
\end{table}

\subsection{Comparison Results}
\label{sec:quad_compare}
\revision{
In this subsection, we compared POPS with other methods in the literature. We first compare to vanilla method and three methods focusing on adaptive sampling according to the initial points:}
\begin{itemize}
    \item \revision{\textit{Vanilla sampling}: Training a model on directly sampled 1500 optimal paths;}
    \item \revision{\textit{Adaptive sampling with large control norms (\textit{AS w. large u})}: A method proposed by \citet{nakamura2021adaptive}, selecting initial points with large gradient norms, equivalent to selecting initial points with large optimal control norms;}
    \item \revision{\textit{Adaptive sampling with large total costs (\textit{AS w. large v})}: selecting initial points whose total costs are large under the latest NN controller;}
    \item \revision{\textit{Adaptive sampling with large value gaps (\textit{AS w. bad v})}: A variant of the SEAGuL algorithm (Sample Efficient Adversarially Guided Learning \citep{landry21seagul}), selecting initial points with large gaps between the learned values and optimal values.}
\end{itemize}

For the three adaptive sampling methods, we sample two points and select one based on the relevant criterion to add an initial point. Each method is run five times with different random seeds.
All methods start with the same initial network as POPS (the policy $\hat{\bm u}_0$ trained on 500 paths) and progressively add 400, 300, and 300 paths, resulting in a final network trained on 1500 paths.

\begin{figure}[t]
\begin{floatrow}
\ffigbox[\FBwidth]{%
  \includegraphics[height =0.24\textwidth]{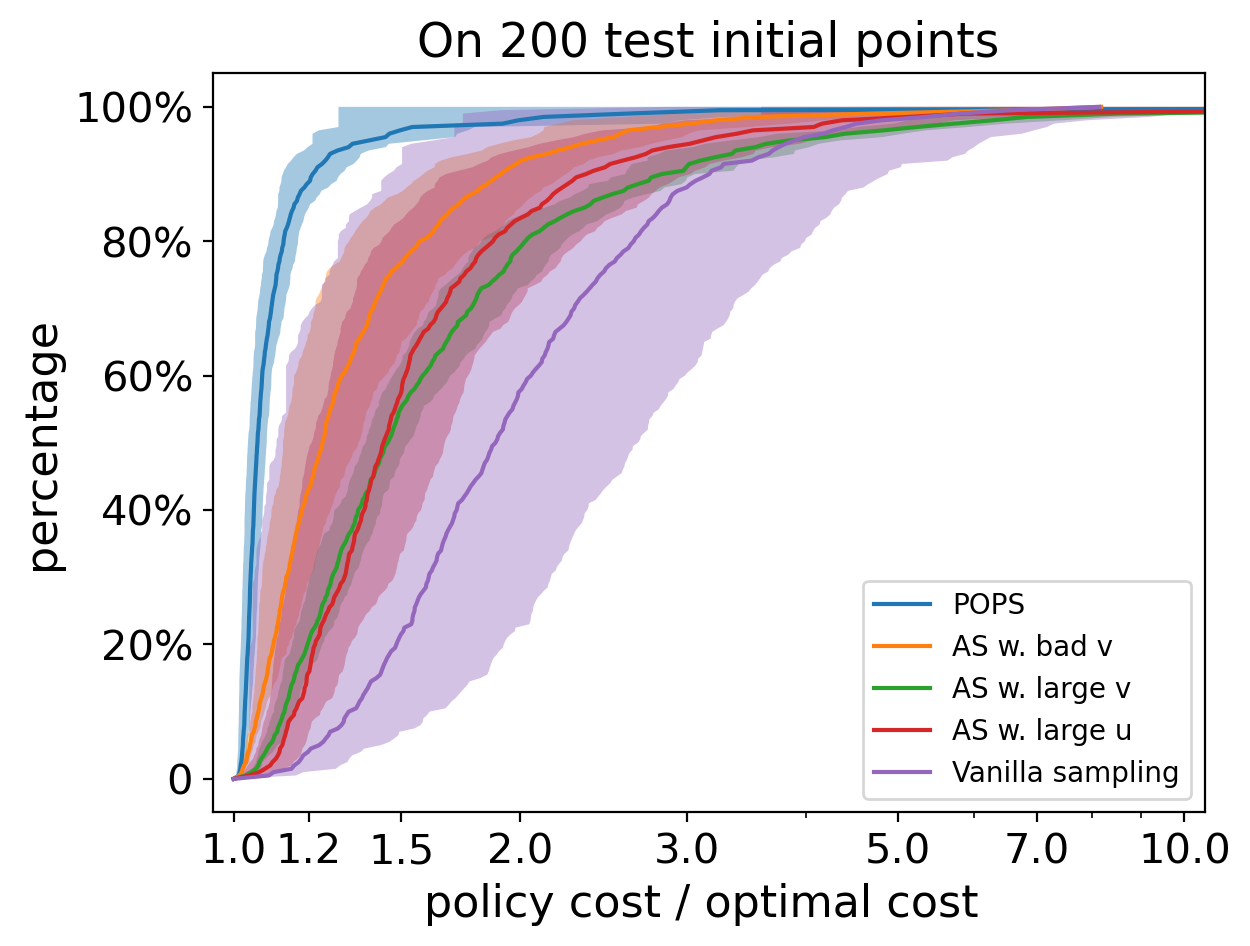}%
}{%
  \caption{Cumulative distribution function of the cost ratio \revision{among 200 test trajectories} across different sampling methods.}\label{fig:quadrotor_cmpfig}%
}
\capbtabbox{%
\footnotesize
\begin{tabular}{lccc}
\toprule
\textbf{Methods}  & \textbf{Mean}  & \textbf{90\%} & \textbf{Median} \\ \midrule
\textbf{POPS}    &  1.17 ($\pm$ 0.09) & 1.22 ($\pm$ 0.09) & 1.06 ($\pm$ 0.02) \\
\textbf{AS w. bad v}&1.45 ($\pm$ 0.18) & 1.93 ($\pm$ 0.33)  & 1.24 ($\pm$ 0.12) \\
\textbf{AS w. large v}   &1.94 ($\pm$ 0.09) & 2.97 ($\pm$ 0.38) & 1.46 ($\pm$ 0.08) \\
\textbf{AS w. large u}  &1.75 ($\pm$ 0.33) & 2.39 ($\pm$ 0.67)  & 1.44 ($\pm$ 0.23) \\
\textbf{Vanilla sampling}   &2.15 ($\pm$ 0.96) & 3.18 ($\pm$ 1.71)  & 1.89 ($\pm$ 0.78) \\
\bottomrule
    \end{tabular}
}{%
  \caption{Statistical results of different sampling methods. \revision{The columns ``Mean'' ``90\%'', and ``Median'' represent the cost ratio statistics across 200 test initial points.}}\label{fig:quadrotor_cmptable}
}
\end{floatrow}
\end{figure}

Data generation is the most time-intensive part of the algorithm, and all methods being compared solve 1500 open-loop problems, except for \textit{AS w. bad v}, which solves 2500. In POPS, the time span of optimal paths shortens with each iteration, reducing computation time compared to other methods, where the time span remains fixed at $T$.
The cumulative distribution functions of cost ratios for these methods, shown in Figure~\ref{fig:quadrotor_cmpfig}, demonstrate the superior performance of POPS. Additional statistics are provided in Table \ref{fig:quadrotor_cmptable}.

\revision{We also compare POPS with \dagger, which augments the training dataset by sampling intermediate states rather than initial points.
In \dagger, we use the same temporal grid points at $t=10$ and $14$ as in POPS. The results, summarized in Table \ref{tab:quadrator_cost_dagger}, show that \dagger performs similarly to POPS when using 500 initial trajectories. However, when the number of trajectories is reduced from 500 to 300, \dagger shows a more significant performance drop.
The main reason is that \dagger demands enough data at the beginning to have a good initial controller capable of exploring the state space over the entire time interval.
However, in complicated control problems, this is often not feasible, making it necessary to employ a more gradual and adaptive sampling strategy, like POPS, to improve the controller effectively.}
\begin{table}[h]
    \centering
    \begin{tabular}{c|cc}
        \toprule
&\multicolumn{2}{c}{\textbf{Number of trajectories}}\\
 \textbf{Methods}&500 &300 \\
\midrule
        \textbf{POPS} &  $1.17$ ($\pm~ 0.09$)& $1.34$ ($\pm~0.13$)  \\
        \textbf{\dagger}  &  $1.19$ ($\pm~0.06$)& $1.50$ ($\pm~0.13$)\\
        \bottomrule
    \end{tabular}
    \caption{\revision{Average cost ratio on 200 test points of controllers trained by POPS and \dagger.
    All models are trained with time grids $0<10<14<16$.
    Cost ratios are clipped at $10.0$ for each test trajectory. For each setting, we repeat 5 times with different random seeds. \dagger performs similarly to POPS with 500 trajectories but shows a noticeable drop in performance when reduced to 300.
    }}
    \label{tab:quadrator_cost_dagger}
\end{table}

\section{The Optimal Reaching Problem of a 7-DoF Manipulator}
\label{sec:reaching_problem}

In this section, we consider the optimal reaching problem on a 7-DoF torque-controlled manipulator, the KUKA LWR iiwa R820 14 \citep{KUKA, kuka-platform}. We formulate the dynamics of the manipulator as 
\[
    \dot{\bm x} = \bm f(\bm x, \bm u) = (\bm v, \bm a(\bm x, \bm u))
,\] 
where $\bm x = (\bm q, \bm v) \in\bR^{14}$, $\bm q\in\bR^7$ is the joint angles, $\bm v=\dot{\bm q}\in\bR^{7}$ is the joint velocities, $\ddot{\bm q} = \bm a(\bm x, \bm u) \in \bR^7 $ is the acceleration of joint angles, and $\bm u\in\bR^7$ is the control torque.
The forward dynamics $\bm a$ is detailed in Appendix \ref{app:dynamics_manipulator}.

\begin{figure}[t]
     \centering
     \includegraphics[width=0.5\textwidth]{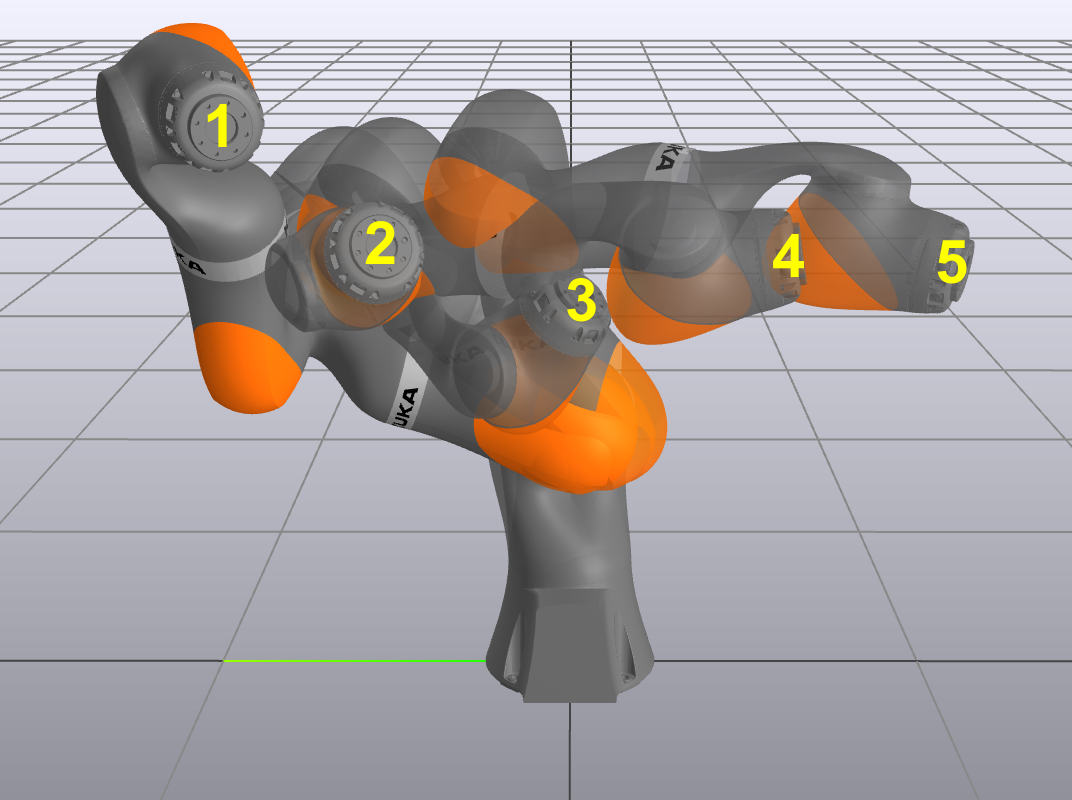}
        \caption{
        An illustration of the reaching problem of the manipulator.
        The solid manipulator demonstrates its initial position.
        We label the end effectors of the five instances of robots by ``1,2,3,4,5'' to indicator the position of the robot at different times  $t_1=0.0, t_2=0.2, t_3=0.4, t_4=0.6, t_5=0.8$.
        }
\label{fig:manipulator_illustration}
\end{figure}

Our goal is to find the optimal torque $\bm u \in \mathcal{U}\subset \bR^{7}$ that drives the manipulator from $ \bm x_0 $ to $\bm x_1 $ in $T=0.8$ seconds and minimizes a quadratic type cost. \revisionjmla{This cost reflects two competing objectives: tracking a desired terminal state and minimizing control-related effort along the way. Specifically, the running cost penalizes both dynamic acceleration and deviation from the gravity-compensating torque $\bm u_1$, while the terminal cost enforces accurate arrival at the target state $\bm x_1$.} See Figure \ref{fig:manipulator_illustration} for an illustration of the task.
In the experiments, we take $\bm x_0 = (\bm q_0, \bm 0), \bm x_1 = (\bm q_1, \bm 0)$ with $\bm q_0 = [1.68,  1.25,  2.44 , -1.27, -0.98, 1.12, -1.36]\transpose$ and $\bm q_1 = [2.77,  0.58,  1.54, -1.70, -2.17,  0.08, -2.58]\transpose$.
The initial positions $\bm q$ are sampled uniformly and independently in a $7$-dimensional cube centered at $\bm{q}_0$ with side length $0.02$. Initial velocities $\bm v$ are set to zero. The running cost and terminal cost are
\begin{align*}
    L(\bm x,\bm u) = \bm a(\bm x,\bm u)\transpose Q_{\bm a} \bm a(\bm x,\bm u) + (\bm u-\bm u_1)\transpose Q_{\bm u} (\bm u-\bm u_1),\,
        M(\bm x) = (\bm x-\bm x_1)\transpose Q_f(\bm x-\bm x_1),
\end{align*}
where  $\bm u_1$ is the torque to balance gravity at state $ \bm x_1 $, \ie $\bm a(\bm x_1, \bm u_1) = \bm{0}$.
Under this setting, $(\bm x_1, \bm u_1)$ is an equilibrium of the system, \ie $\bm f(\bm x_1, \bm u_1) = (\bm v_1, a(\bm x_1, \bm u_1)) = \bm{0}$.  We take $Q_{\bm a}= 0.005 I_7, Q_{\bm u} = 0.025 I_7, Q_f=25000 I_{14} $ where we use large weights $Q_f$ to ensure the reaching goal is approximately achieved.

To obtain training data,  we solve the open-loop control problem use differential dynamic programming \citep{DDP} implemented in the Crocoddy library \citep{mastalli20crocoddyl}.
\revisionjmla{
Instead of directly applying the open-loop solver to each collected initial state, we first sample a mini-batch of initial states around $\bm q_0$ from the same distribution and solve the open-loop problem for each. We then select the solution with the lowest cost and use it as an initial guess to generate all remaining trajectories. This strategy accelerates data generation, helps avoid poor local minima, and leads to a $100\%$ success rate in solving the open-loop control problem.
}
In the simulation and open-loop solver, we take time step $\Delta t=0.001$ and use the semi-implicit Euler discretization.
Each trajectory has $T/\Delta t = 800$ data points that are pairs of $15$-dimensional input states (including time) and $7$-dimensional output controls.

The backbone network for this example is the QRNet \citep{qrnet1,nakamura2021neural}. %
QRNet exploits the solution corresponding to the linear quadratic regulator (LQR) problem at equilibrium and thus improves the network performance around the equilibrium~\citep{hu2023learning}.
The usage of a different network structure in this example also demonstrates the genericness/versatility of POPS.
We leave the details for QRNet to Appendix \ref{sec:QRnet}. 

All the QRNets $\bm u^{\QR}$ are trained  with the Adam optimizer \citep{adam} with learning rate $0.001$, batch size $256$ and epochs $2000$.
We utilize a fully-connected network with $6$ hidden layers; each layer has $128$ neurons.
\revision{Compared to the quadrotor problem, here we used a deeper network and replaced the bounded \texttt{tanh} activation with the unbounded \texttt{ELU} \citep{elu} in the last three layers to handle the manipulator problem’s increased complexity and significantly varying scale.}
During iterations, all networks are trained from scratch, \ie a new network with random weights instead of loading weights from the previous iteration.

We evaluate networks trained in six different ways: four using Algorithm \ref{alg:main} with different temporal grid points (\textit{POPS1}--\textit{POPS4}), and two using the vanilla sampling method with $300$ (\textit{Vanilla300}) and $900$ (\textit{Vanilla900}) trajectories.
\textit{POPS1}--\textit{POPS4} are trained with an initial training data of $100$ trajectories and undergo three iterations ($K=3$), \ie, each of them requires solving the open-loop solution $300$ times in total.

\begin{figure}[t]
\centering
\begin{subfigure}[b]{0.32\textwidth}
    \centering\includegraphics[width=\textwidth]{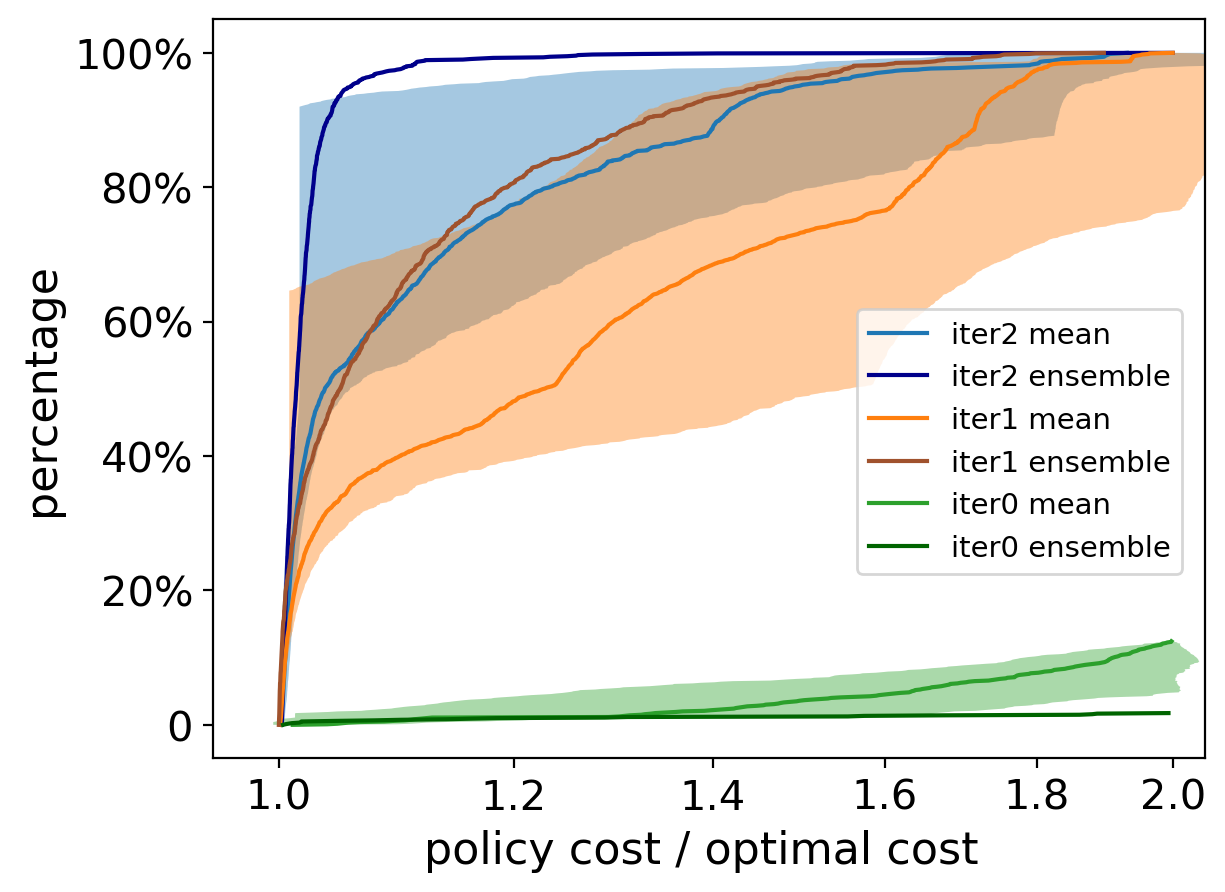}
\end{subfigure}
\begin{subfigure}[b]{0.32\textwidth}
    \centering\includegraphics[width=\textwidth]{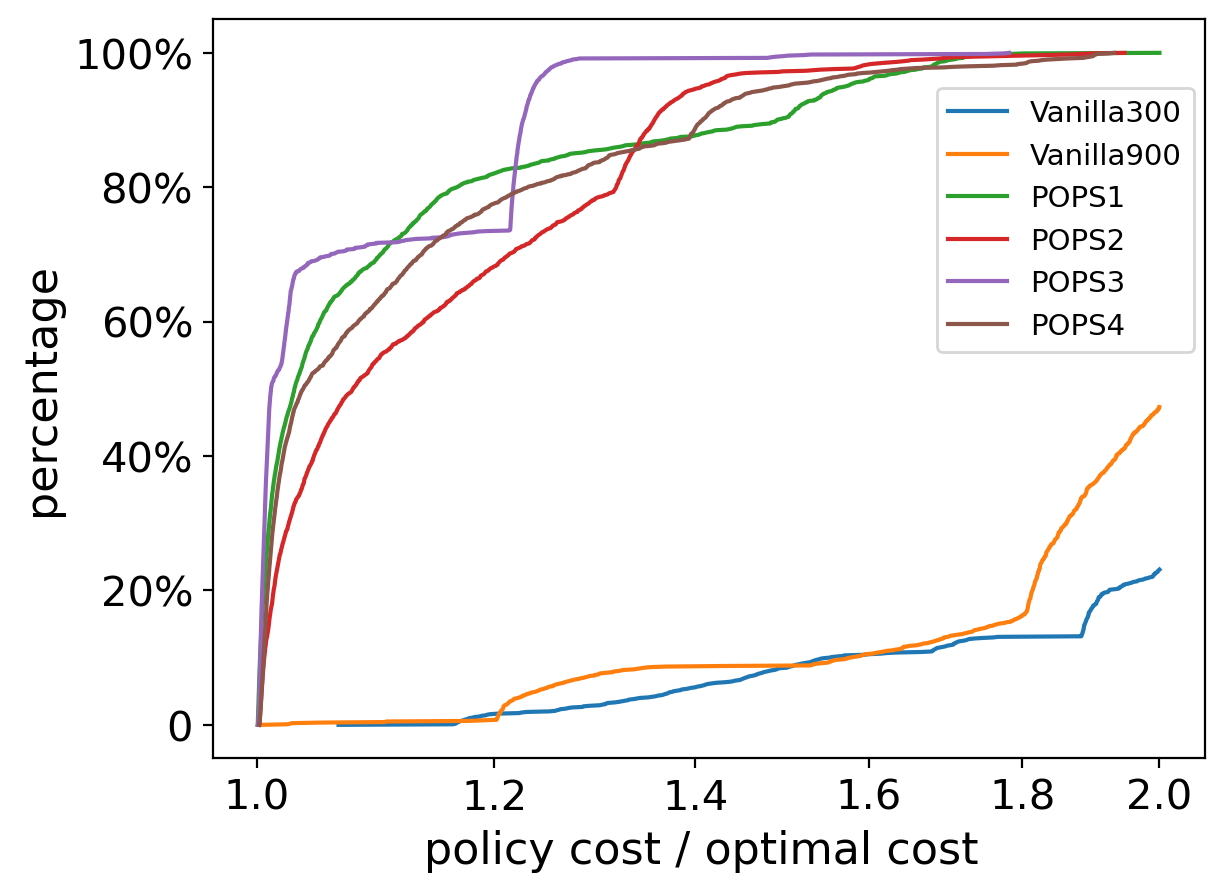}
\end{subfigure}
\begin{subfigure}[b]{0.32\textwidth}  
    \centering \includegraphics[width=\textwidth]{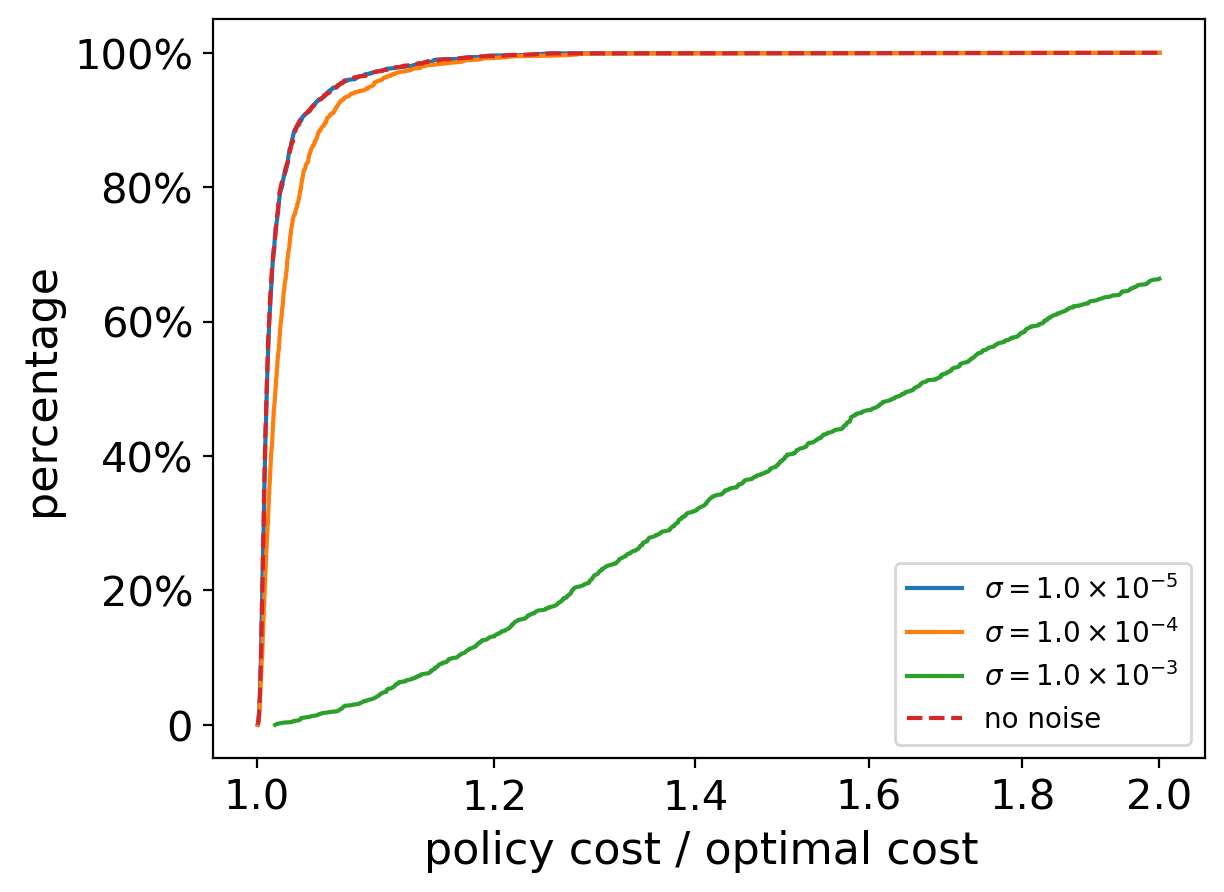}
\end{subfigure}
\caption
{Cumulative distribution functions of cost ratios (\idealRatio)
under different iterations when training \textit{POPS4} (left), different training schemes (middle) and different intensities of measurement noises (right).
The shaded regions in the left plot represent the mean $\pm$ standard deviation. \revision{The percentages are computed based on the cost ratios of 1200 test trajectories}.
}
\label{fig:manipulator}
\end{figure}

Each experiment repeat five times independently and evaluate them on the test dataset comprising 1200 trajectories.
We again plot the cumulative distribution functions of cost ratios (clipped at 2.0) between the NN-controlled cost and optimal cost in Figure \ref{fig:manipulator}.
As reported in Table \ref{tab:reaching_costs}, we find that adding more data in the vanilla sampling method has very limited effects on improvement while POPS greatly improves the performance. Again, such improvement is robust to different choices of temporal grid points (\textit{POPS1}--\textit{POPS4}).
Besides, we also try augmenting the dataset with newly collected data instead of replacing them, as detailed in Section~\ref{sec:conclusion}.
Through the comparison between \textit{POPS3} and \textit{POPS3*}, we find that the alternative approach does not bring further improvement.

We additionally evaluate the NN controllers in the presence of measurement errors by adding a noise uniformly sampled from $[-\sigma, \sigma]^{14}$ to the state $\bm x$ of the network input.
See Figure \ref{fig:manipulator} (right) for the results on the best model trained from \textit{POPS1}--\textit{POPS4}.
The NN controller performs well at $\sigma=10^{-4}$, and there are more than $60\%$ of cases in which our controller achieves a ratio less than $2.0$ at $\sigma=10^{-3}$.

\begin{table}[ht]
\centering
\begin{tabular}{c|c|ccc}
\toprule
&\textbf{time grid} & iter 0&iter 1 &iter 2\\
\midrule
\textbf{Vanilla300} &&1.92 ($\pm$ 0.06)&&\\
\textbf{Vanilla900} &&1.89 ($\pm$ 0.09)&&\\
\textbf{POPS1} & 0.16 - 0.48 - 0.8 &1.81 ($\pm$ 0.16) &  1.29 ($\pm$ 0.21) &  1.09 ($\pm$ 0.08)\\
\textbf{POPS2}  & 0.16 - 0.56 - 0.8 &1.94 ($\pm$ 0.07) &  1.34 ($\pm$ 0.17) &  1.15 ($\pm$ 0.14)\\
\textbf{POPS3} & 0.16 - 0.64 - 0.8 & 1.92 ($\pm$ 0.08) &  1.37 ($\pm$ 0.21) &  1.07 ($\pm$ 0.11)\\
\textbf{POPS3*} & 0.16 - 0.64 - 0.8 & 1.96 ($\pm$ 0.05) &  1.40 ($\pm$ 0.17) &  1.24 ($\pm$ 0.28)\\
\textbf{POPS4} &0.16 - 0.72 - 0.8 &1.96 ($\pm$ 0.03) &  1.29 ($\pm$ 0.18) &  1.13 ($\pm$ 0.10) \\
\bottomrule
\end{tabular}
\caption{The mean ratio of policy costs / optimal costs of the optimal reaching problem of the manipulator.
The ratio has been clipped at $2.0$ for each test trajectory. The vanilla300/900 correspond to networks trained on 300/900 optimal trajectories, respectively.
The choices of temporal grid points for adaptive sampling in the remaining rows can be inferred by the location of columns. For example, POPS1 has temporal grid points $0<0.16<0.48<0.8$.
POPS3* has the same temporal grid points as POPS3 except that it augments the dataset directly instead of replacing them, as discussed in Section~\ref{sec:conclusion}.}
\label{tab:reaching_costs}
\end{table}

\subsection{Comparison with \dagger}\label{sec:reaching_compare}
\revision{
In this subsection, we compare our method with \dagger. In \dagger, we use the same temporal grid points at $t=0.16$ and $t=0.64$ as \textit{POPS3}.
The \dagger algorithm achieves a policy cost / optimal cost ratio of $1.049$ on the test dataset, which is close to that achieved by POPS.}

\revision{
We then increase the difficulty of the control problem by enlarging the moving distance. 
Specifically, we change the center of the initial position and the terminal position to 
\begin{align*}
    &\bm q_0 = [1.60, 1.30, 2.70, -0.85, -1.90, 0.95, -1.60]\transpose, \\
    \text{and~~} &\bm q_1 =[2.75, 0.60, 2.00, -1.55, -2.15, 0.00, -2.60]\transpose,
\end{align*}
respectively.
We also increase the size of the initial dataset from 100 trajectories to 200 trajectories. Each network is trained with 1500 epochs. The other settings remain unchanged.}

\begin{figure*}[t]
\centering
\begin{subfigure}[b]{0.45\textwidth}
    \centering\includegraphics[width=\textwidth]{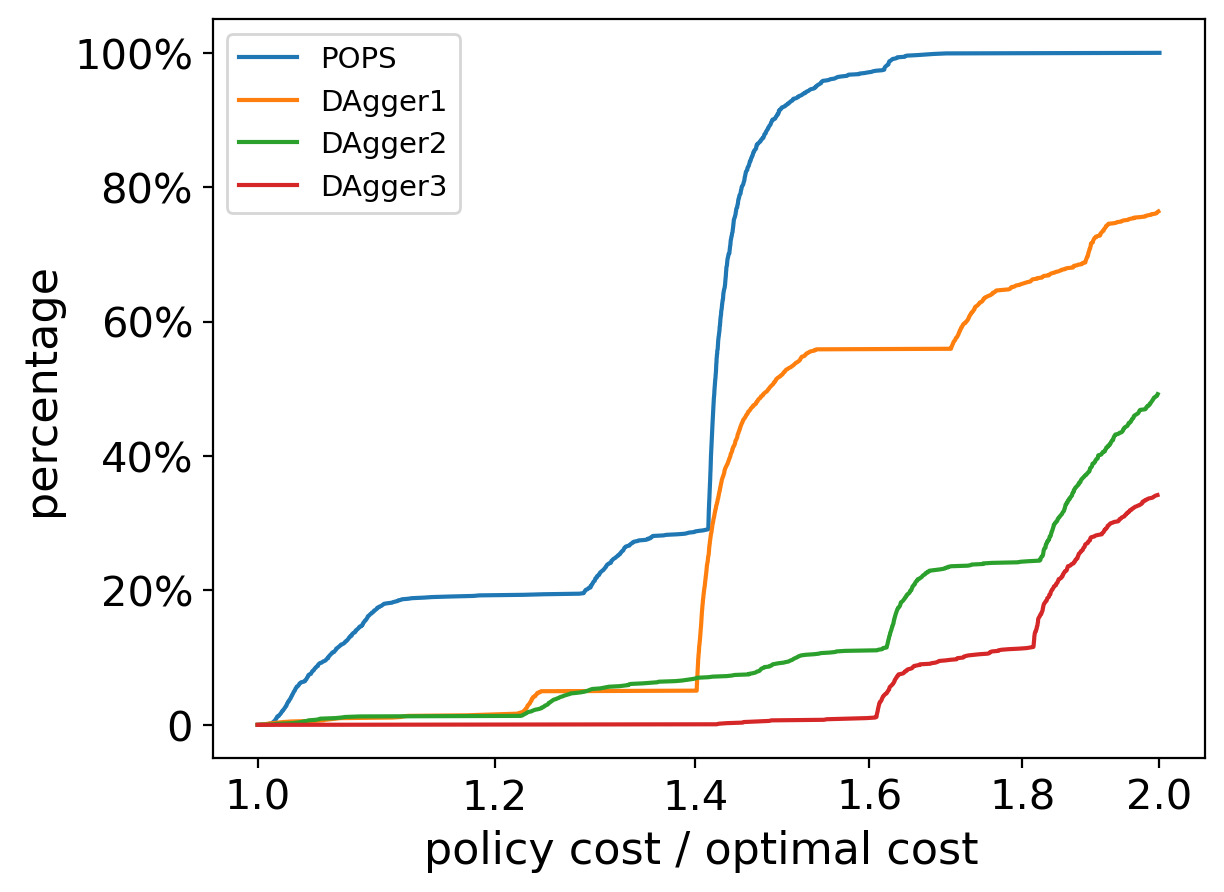}
\end{subfigure}
\quad
\begin{subfigure}[b]{0.45\textwidth}
    \centering\includegraphics[width=\textwidth]{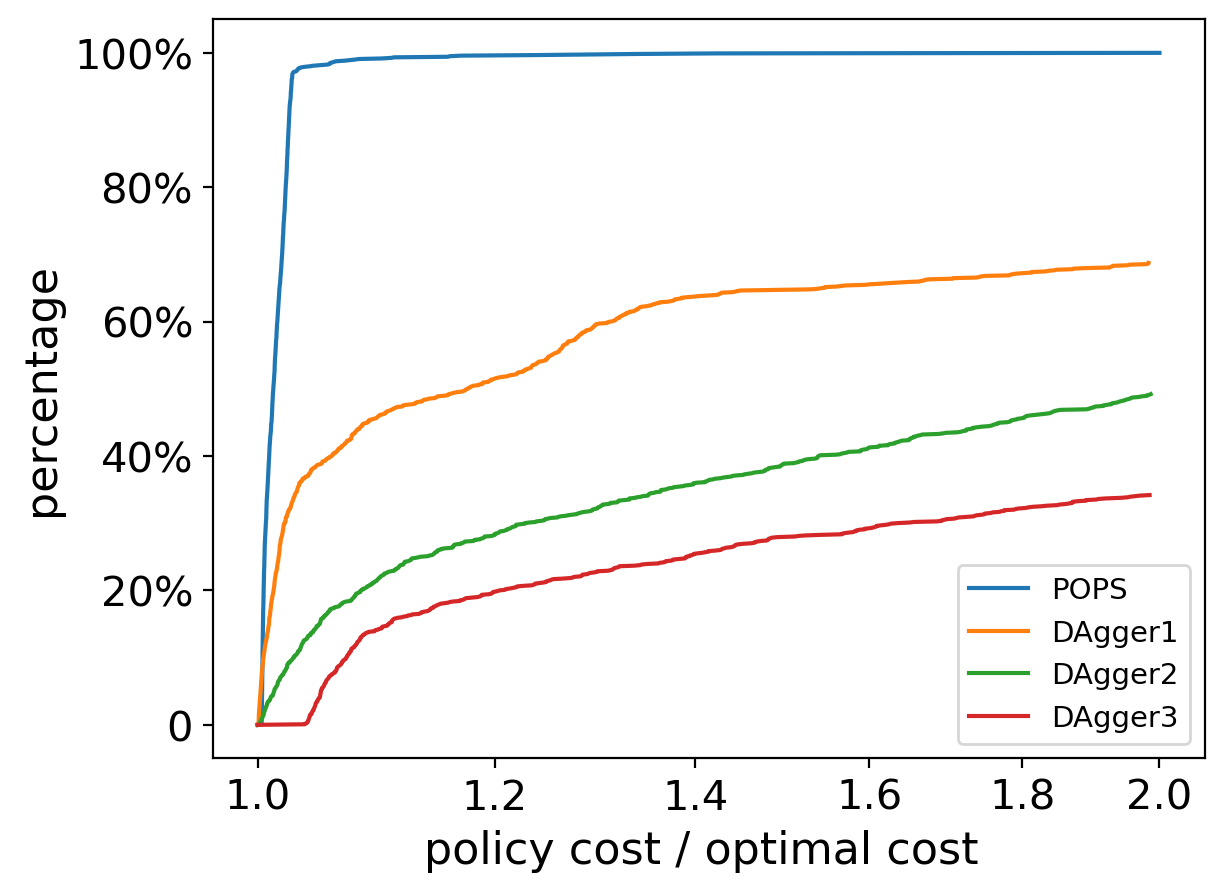}
\end{subfigure}
\caption{Cumulative distribution functions of average cost ratio over 5 independent experiments (left) and the cost ratio of the best controller among 5 independent experiments (right)
under the proposed method and \dagger for the optimal reaching problem with a larger moving distance. \revision{The percentages are computed based on the cost ratios of 1200 test trajectories}.
}
\label{fig:manipulator_ratio_difficult}
\end{figure*}

\revision{
For a comprehensive comparison, we choose two different temporal grids points. \textit{\dagger1} uses the same temporal grids as \textit{POPS1} while \textit{\dagger2} uses the same temporal grid points as \textit{POPS3}.
Each experiment has been run 5 times independently, and we report their average and best performance.
The results are summarized in Table \ref{tab:reaching_difficult_costs} and Figure \ref{fig:manipulator_ratio_difficult}.
As we can see, POPS is capable of finding a closed-loop controller with an average ratio between policy cost and optimal cost achieving $1.0155$.
However, the \dagger algorithm cannot yield such a satisfactory result.}

\revision{
As discussed in Section \ref{sec:quad_compare}, we maintain that \dagger requires a strong initial controller for effective state exploration; otherwise, states sampled by a poor controller might even deteriorate the performance.
The results in Table \ref{tab:reaching_difficult_costs} support this view.
First, we observe that \textit{\dagger1} performs similarly to the network trained in the second iteration of POPS, suggesting that additional data sampled at the later time $t=0.48$ offers little benefit. Furthermore, comparing \textit{\dagger1} and \textit{\dagger2}, where \textit{\dagger1} has an earlier final grid point ($0.48$ vs $0.64$) and better performance (mean ratio $1.8528$ vs $1.6327$), reveals that late-time sampling actually worsens the performance. Lastly, we performed an additional iteration of \dagger (\textit{\dagger3}), which requires adding 400 more trajectories. However, this led to worse performance, as the new states were sampled by a less effective controller.}

\begin{table}[ht]
    \centering
    \begin{tabular}{c|c|ccc}
        \toprule
&&\multicolumn{3}{c}{\textbf{Temporal grid points ($T=0.8$)}}\\
&\textbf{$\#$ of iterations} & 0.16 & 0.48 &  0.64 \\
\midrule
        \textbf{POPS} & 3  & 1.6592 (1.5417) && 1.3610 (1.0155) \\
        \textbf{\dagger1} & 2 & & 1.6327 (1.4004) &   \\
        \textbf{\dagger2} & 2 &  && 1.8528 (1.6364)   \\
        \textbf{\dagger3} & 3 &  && 1.9293 (1.7478) \\
        \bottomrule
    \end{tabular}
    \caption{
    The mean ratio between policy costs and optimal costs of the reaching problem with a larger moving distance. In each cell, the first number is averaged over 5 independent experiments and the number in parenthesis is the average ratio achieved by the best controller among them.
The ratio has been clipped at $2.0$ for each test trajectory. 
Both \textit{POPS} and \textit{\dagger2} have temporal grid points $0<0.16<0.64<0.8$.
\textit{\dagger3} repeats \textit{\dagger2} for one more \dagger iteration on the same temporal grid points.
    }
    \label{tab:reaching_difficult_costs}
\end{table}

\section{Conclusion and Future Work}
\label{sec:conclusion}
In this work, we propose the progressive optimal path sampling method to overcome the distribution mismatch problem in the supervised-learning-based approaches for the closed-loop optimal control problem. Both theoretical and numerical results show that the POPS method significantly improves the performance of the learned NN controller and outperforms other adaptive sampling methods. 

There are several directions worth exploring in future work. \revisionjmlb{One design choice in POPS is the selection of temporal grid points for adaptive sampling. In this paper, we use a heuristic that selects a common grid point for all trajectories, based on when the distance between the training data and NN-controlled trajectories begins to increase rapidly (see Figure \ref{fig:sim2train} for an example). We observe that POPS performs well with this strategy. A recent work \cite{hu2023learning} proposes a more systematic approach that selects grid points individually for each trajectory, based on the first significant deviation from its optimal counterpart. It would be worthwhile to further develop and generalize such strategies.
}%

Another direction is to design more effective approaches utilizing the training data. In Algorithm~\ref{alg:main} (lines 8--9), at each iteration, we replace parts of the training data with the newly collected data, and hence some optimal labels are thrown away, which are costly to obtain. An alternative choice is to augment data directly, \ie, setting $S_i = \hat{S}_i \bigcup S_{i-1}$ in line 9. Numerically, we observe that this choice gives similar performance to the version used in Algorithm~\ref{alg:main}, which suggests that so far the dropped data provides little value for training. However, it is still possible to find smarter ways to utilize them to improve performance.
\revisionjmlb{Beyond data efficiency, it is important to evaluate POPS on more challenging classes of optimal control problems, such as those involving state or control constraints. Another important challenge is the presence of multi-modal optimal policies. One strategy to address this is to pre-train a neural network and use it to warm-start an open-loop optimal-control solver, thereby restricting the training set to a single mode that the network can accurately approximate \cite{nakamura2021adaptive}. A complementary direction is to model the policy stochastically with a generative approach, such as a diffusion model, so that multiple modes are captured implicitly \cite{janner2022diffuser,chi2023diffusion,hegde2023generating,domingo2024stochastic}. It is of interest to explore how POPS can be integrated with neural-network warm starts and stochastic policy representations.}

Finally, POPS can be straightforwardly applied to learning general dynamics from multiple trajectories as the controlled system under the optimal policy can be viewed as a special dynamical system. It is of interest to investigate its performance in such general settings. Theoretical analysis beyond the LQR setting is also an interesting and important problem.

\appendix

\section{Detailed Setting and Proof in Section \ref{sec:toy_example}}\label{appendix_lqr}
We first present the settings for the vanilla supervised-learning-based method, \dagger, and POPS when applied to the LQR problem in Section \eqref{sec:toy_example}.

\paragraph{Vanilla method.} For the vanilla method, we first randomly sample $NT$ initial states $\{\tilde{x}_j^v\}_{j=1}^{NT}$\footnote{Through the LQR analysis, all symbols having a hat are open-loop optimal paths sampled for training, \eg $\hat u_j^v, \hat x_j^v, \hat u_{i,j}^p, \hat x_{i,j}^p$.
Let $\tilde{x}$ denote a single state instead of a state trajectory.
The clean symbol $x$ without hat or tilde is the IVP solution generated by specific controllers which are specified in the subscript; \eg $x_o,x_v,x_d,x_p$ are trajectories generated by $u_o, u_v, u_d, u_p$ which are optimal, vanilla, \dagger and POPS controllers, respectively. 
The positive integer $j$ in the subscript always denotes the index of the optimal path.
Symbols with superscript $i$ are related to the $i$-th temporal grid points in \dagger or POPS.} from a standard normal distribution where $N$ is a positive integer (recalling $T$ is a positive integer). Then $NT$ approximated optimal paths are collected starting at $t_0=0$:
\begin{equation}\label{vanilla_data}
    \hat{u}_j^v(t) = -\frac{T}{T^2 + 1}\tilde x_j^v + \epsilon Z_j^v,\quad \hat{x}_j^v(t) = \frac{T(T-t)+1}{T^2 + 1}\tilde x_j^v + \epsilon t Z_j^v,
\end{equation}
where $\{Z_j^v\}_{j=1}^{NT}$ are \textit{i.i.d.} normal random variables with mean $m$ and variance $\sigma^2$, and independent of initial states.

Finally, the parameters $\theta$ are learned by solving the following least square problems:
\begin{equation}\label{LS_vanilla}
    \min_{\theta}\int_{0}^{T}\sum_{j=1}^{NT}|\hat{u}_j^v(t) - u_\theta(t,\hat{x}_j^v(t))|^2\rmd t.
\end{equation}
For the first and second models, we optimize $\theta_t$ independently for each $t$:
\begin{equation*}
   \theta_t = \argmin_{ b}\sum_{j=1}^{NT}|\hat{u}_j^v(t) + \frac{T}{T(T-t)+1}\hat{x}_j^v(t) - b|^2 \,\text{ or }\,\theta_t = \argmin_{ (a,b)}\sum_{j=1}^{NT}|\hat{u}_j^v(t) - a\hat{x}_j^v(t) - b|^2.
\end{equation*}
 We will use $u_v$ to denote the closed-loop controller determined in this way. Notice that $\{\hat{x}_j^v(t)\}_{j = 1}^N$ share the same distribution, we have that $\hat{x}^v(t) $ has the same distribution of $\hat{x}_1^v(t)$

\paragraph{\dagger.} In \dagger, we again choose $K=T$ and the temporal grid points $t_i = i$ for $0 \le i \le K$. We first sample $N$ initial points $\{\tilde{x}_{0,j}^d\}_{j=1}^N$ from the normal standard distribution and then generated $N$ approximated optimal paths starting at $t_0 = 0$:
\begin{equation*}
    \hat{u}_{0,j}^d(t) = -\frac{T}{T^2+1}\tilde{x}_{0,j}^d + \epsilon Z_{0,j}^d, \quad \hat{x}_{0,j}^d(t) = \frac{T(T-t)+1}{T^2+1}\tilde{x}_{0,j}^d + \epsilon t Z_{0,j}^d,
\end{equation*}
where $\{Z_{0,j}^d\}_{j=1}^N$ are \textit{i.i.d.} normal random variables and independent with initial states whose mean is $m$ and variance is $\sigma^2$. We then train the closed-loop controller $u_0$ by solving the following least square problems:
\begin{equation*}
     \min_{\theta}\int_{0}^{T}\sum_{j=1}^{N}|\hat{u}_{0,j}^d(t) - u_0(t,\hat{x}_{0,j}^d(t))|^2\rmd t.
\end{equation*}

Then, we use $u_d^0$ to solve the IVPs on the whole time horizon $[0,T]$ with initial states $\{\tilde{x}_{0,j}^d\}$:
\begin{equation}\label{ivp_DAGGER}
    \dot{x}_{0,j}^d(t) = u_0(t,x_{0,j}^d(t)),\quad x_{0,j}^d(0) = \tilde{x}_{0,j}^d,\, 1\le j \le N,
\end{equation}
and collect $\{\tilde{x}_{i,j}^d\}_{j=1}^{N}$ as $\tilde{x}_{i,j}^d:=x_{0,j}^d(i)$ for $i = 1,2,\dots,T-1$. At each time step $t_i = i$, we then compute $N$ approximated optimal paths starting from $\{\tilde{x}_{i,j}^d\}_{j=1}^N$:
\begin{equation}\label{DAGGER_data}
    \hat{u}_{i,j}^d(t) = -\frac{T}{T(T-i)+1}\tilde{x}_{i,j}^d + \epsilon Z_{i,j}^d, \quad \hat{x}_{i,j}^d(t) = \frac{T(T-t)+1}{T(T-i)+1}\tilde{x}_{i,j}^d + (t-i)\epsilon Z_{i,j}^d, \,t \in [i,T],
\end{equation}
where $\{Z_{i,j}^d\}_{1\le i \le T-1}$ are \textit{i.i.d.} normal random variables and independent with $\{\tilde{x}_{0,j}^d\}_{j=1}^{N}$ and $\{Z_{0,j}^d \}_{j=1}^N$ whose mean is $m$ and variance is $\sigma^2$. Finally, we collect the optimal paths $\{(\hat{u}_
{i,j}^d,\hat{x}_{i,j}^d)\}_{0\le i \le T-1, 1\le j \le N}$ to train the closed-loop controller $u_d$ by solving the following least square problems:
\begin{equation}\label{least_square_DAGGER}
    \min_{\theta} \int_{i}^{i+1}\sum_{k=0}^{i}\sum_{j=1}^{N}|\hat{u}_{k,j}^d(t) - u_\theta(t,\hat{x}_{k,j}^d(t))|^2\rmd t
\end{equation}
for $i = 0,1,\dots, T-1$. In the theoretical part, we will only do a single iteration and use $u_d$ to denote the policy learned by \eqref{least_square_DAGGER}. In the numerical part, we can replace $u_0$ by the policy learned by \eqref{least_square_DAGGER} and repeat this process multiple times.

\paragraph{POPS.} In POPS, we choose $K = T$ and the temporal grid points $t_i = i$ for $0 \le i \le K$.
We first sample $N$ initial points $\{\tilde x_{0,j}^p\}_{j=1}^{N}$ from the normal standard distribution, denote the parameters optimized at $i$-th iteration as $\theta^{i}$ and initialize $\theta^{-1} = 0$. At the $i$-th iteration ($0 \le i \le T-1$), we use $ u_{\theta^{i-1}}$ to solve the IVPs on the time horizon $[0,i]$
\begin{equation}\label{IVP_iteration}
    \dot{x}_{i,j}^p(t) = u_{\theta^{i-1}}(t,x_{i,j}^p(t)), \quad x_{i,j}^p(0) = \tilde x_{0,j}^p,\, 1\le j \le N,
\end{equation}
and collect $\{\tilde x_{i,j}^p\}_{j=1}^N$ as $\tilde x_{i,j}^p \coloneqq x_{i,j}^p(i)$. 
We then compute $N$ approximated optimal paths starting from $\{\tilde x_{i,j}^p\}_{j=1}^{N}$ at $t_i = i$:
\begin{equation}\label{adaptive_data}
    \hat{u}_{i,j}^p(t) = -\frac{T}{T(T-i) + 1}\tilde x_{i,j}^p + \epsilon Z_{i,j}^p,\quad \hat{x}_{i,j}^p(t) = \frac{T(T-t)+1}{T(T-i) + 1}\tilde x_{i,j}^p + (t-i)\epsilon  Z_{i,j}^p, \, t \in[i,T]
\end{equation}
where $
\{Z_{i,j}^p\}_{0\le i\le T-1,1\le j \le N}$ are \textit{i.i.d.} normal random variables with mean $m$ and variance $\sigma^2$, and independent of $\{\tilde x_{0,j}^p\}_{j=1}^{N}$.
Note that $\hat u_{i,j}^p$ and $\hat x_{i,j}^p$ are only defined in $t\in[i, T]$ (for $i\ge 1$), we then fill their values in interval $[0, i)$ with values from previous iteration,
\begin{equation}\label{adaptive_iteration}
    \hat{u}_{i,j}^p(t) = \hat{u}_{i-1,j}^p(t),\quad \hat{x}_{i,j}^p(t) = \hat{x}_{i-1,j}^p(t),\quad t \in [0,i).
\end{equation}
Finally, we solve the least squares problems to determine $\theta^i$:
\begin{equation}
\label{LS_adaptive}
    \min_{\theta}\int_{0}^T \sum_{j=1}^{N}|\hat{u}_{i,j}^p(t) - u_{\theta}(t,\hat{x}_{i,j}^p(t))|^2\rmd t.
\end{equation}
We will then use $u_p$ to denote $u_{\theta^{T-1}}$, the policy learned in the last iteration. Again, notice that $\hat{x}_{T-1,j}^p(t)$ share the same distribution, we know that $\hat{x}^p(t)$
has the same distribution of $\hat{x}_{T-1,1}^p(t)$.

Below, we give the proof of Theorem \ref{thm_1}:
\subsection{Vanilla Sampling Method}
We first give the closed-form expressions of $u_v$ using Model 1~\eqref{lqr_model1}. Recalling $\hat{u}_j^v(t)$ and $\hat{x}_j^v(t)$ given in equation \eqref{vanilla_data}, we have
\begin{equation*}
    \hat{u}_j^v(t) = -\frac{T}{T(T-t)+1}\hat{x}_j^v(t) +\frac{T^2+1}{T(T-t)+1}\epsilon Z_j^v, \, 1\le j \le NT.
\end{equation*}
Therefore, recalling $u_v$ is learned through the least square problem \eqref{LS_vanilla}, we have
\begin{equation}\label{eq:vanilla_control}
    u_v(t,x) = -\frac{T}{T(T-t)+1}x + \frac{T^2+1}{T(T-t)+1}\epsilon \bar{Z}^v,
\end{equation}
where
\begin{equation*}
    \bar{Z}^v = \frac{1}{NT}\sum_{j=1}^{NT}Z_j^v.
\end{equation*}
\paragraph*{Distribution Difference.} Using the control~\eqref{eq:vanilla_control}, we have
\begin{equation*}
    \dot{x}_v(t) = -\frac{T}{T(T-t)+1}x_v(t) + \frac{T^2+1}{T(T-t)+1}\epsilon \bar{Z}^v,\; x_v(0) = \txinit.
\end{equation*}
Solving this ODE gives
\begin{equation}\label{IVP_v}
    x_v(t) = \frac{T(T-t)+1}{T^2+1}\txinit + \epsilon t \bar{Z}^v.
\end{equation}

Combining the last equation with the fact that
\begin{equation*}
    \hat{x}_j^v(t) = \frac{T(T-t)+1}{T^2+1}\tilde x_j^v + \epsilon t Z_j^v,
\end{equation*}
$\tilde x_j^v$, $\txinit$ and $\{Z_j^v\}_{j=1}^{NT}$ are independent normal random variables and
\begin{equation*}
    \tilde x_j^v, \txinit \sim \mathcal{N}(0,1), Z_j^v \sim \mathcal{N}(m,\sigma^2),
\end{equation*}
we know that $\hat{x}_j^v(t)$ and $x_v(t)$ are normal random variables with
\begin{equation*}
\mathbb{E} \hat{x}_1^v(t) = \mathbb{E} x_v(t), |\Var|\hat{x}_1^v(t)|^2 - \Var|x_v(t)|^2 | =\sigma^2 (1-\frac{1}{NT})\epsilon^2 t^2.
\end{equation*}

\paragraph*{Performance Difference. } First, with the optimal solution
\begin{equation*}
    u_o(t) = -\frac{T}{T^2+1}\txinit, \quad x_o(T) = \frac{1}{T^2+1}\txinit,
\end{equation*}
we have
\begin{equation*}
    J_o = \frac{1}{T}\int_{0}^{T}\left|\frac{T}{T^2+1}\txinit\right|^2\rmd t + \left|\frac{1}{T^2+1}\txinit\right|^2= \frac{1}{T^2+1}|\txinit|^2.
\end{equation*}
Recalling equation \eqref{IVP_v} and plugging \eqref{IVP_v} into \eqref{eq:vanilla_control}, we know that
\begin{align*}
    &x_v(T) = \frac{1}{T^2+1}\txinit + \epsilon T \bar{Z}^v,\\
    &u_v(t) = -\frac{T}{T^2+1}\txinit + \epsilon \bar{Z}^v.
\end{align*}
Hence,
\begin{equation*}
    J_v = \epsilon^2 |\bar{Z}^v|^2 -\frac{2T}{T^2+1}\txinit\epsilon \bar{Z}^v+ \epsilon^2 T^2 |\bar{Z}^v|^2 + \frac{2T}{T^2+1}\txinit \epsilon \bar{Z}^v + \frac{1}{T^2+1}|\txinit|^2,
\end{equation*}
which gives

\begin{equation*}
    \mathbb{E} J_v  - J_o = 
    (T^2+1)(m^2+\frac{\sigma^2}{NT})\epsilon^2.
\end{equation*}

\subsection{POPS}
Similarly, we first compute $u_p$. Recalling equations \eqref{adaptive_data} and \eqref{adaptive_iteration}, when $0 \le i \le T-1$, $1 \le j \le N$ and $t \in [i,i+1)$, we have
\begin{equation}\label{eq:adaptive_last_iteration}
\begin{aligned}
        &\hat{u}_{T-1,j}^p(t) = \hat{u}_{i,j}^p(t) = -\frac{T}{T(T-i)+1}\tilde x_{i,j}^p + \epsilon Z_{i,j}^p,\\ 
    &\hat{x}_{T-1,j}^p(t) = \hat{x}_{i,j}^p(t) = \frac{T(T-t)+1}{T(T-i)+1}\tilde x_{i,j}^p + (t-i)\epsilon Z_{i,j}^p.
\end{aligned}
\end{equation}
Therefore,
\begin{equation*}
    \hat{u}_{T-1,j}^p(t) = -\frac{T}{T(T-t)+1}\hat{x}_{T-1,j}^p(t) + \frac{T(T-i)+1}{T(T-t)+1}\epsilon Z_{i,j}^p.
\end{equation*}
Hence, recalling $u_p$ is learned through the least square problem \eqref{LS_adaptive}, we have, when $t \in [i,i+1)$
\begin{equation}\label{hat_u_v}
    u_p(t,x) = -\frac{T}{T(T-t)+1}x + \frac{T(T-i)+1}{T(T-t)+1}\epsilon \bar{Z}_i^p,
\end{equation}
where
\begin{equation*}
    \bar{Z}_i^p = \frac{1}{N}\sum_{j=1}^{N} Z_{i,j}^p, \, 0 \le i \le T-1.
\end{equation*}
Equation \eqref{hat_u_v} also holds when $i = T-1$ and $t = T$. 

We then compute the starting points $\{\tilde x_{i,j}^p\}_{0\le i \le T-1, 1 \le j \le N}$ in the POPS method.
By equation \eqref{adaptive_iteration}, we know that when $1\le i \le i' \le T-1$ and $0 \le t < t_i$, $\theta_t^i = \theta_t^{i'}$.
Together with equation \eqref{IVP_iteration}, we know that when $0 \le i \le T-2$, $u_{\theta^i}(t,x) = u_{\theta^{T-1}}(t,x) =  u_p(t,x)$ for $t \in [i,i+1)$, and $x_{i+1,j}^p(t) \equiv x_{i,j}^p(t)$ for $t \in [0,i]$, which implies $x_{i+1,j}^p(i)=x_{i,j}^p(i)=\tilde{x}_{i,j}^p$. Therefore, for $1 \le j \le N$, when $0\le i \le T-2$, we have
\begin{align*}
\begin{cases}
    \dot{x}_{i+1,j}^p(t)&= \displaystyle{u_{\theta^{i}}(t, x_{i+1,j}^p(t)) = u_{p}(t,x_{i+1,j}^p(t))} \\&= \displaystyle{-\frac{T}{T(T-t)+1}x_{i+1,j}^p(t) + \frac{T(T-i)+1}{T(T-t)+1}\epsilon \bar{Z}_i^p},\;  t \in [i,i+1],\\
    x_{i+1,j}^p(i)&= \tilde x_{i,j}^p.
\end{cases}
\end{align*}

Solving the above ODE, we get the solution
\begin{align*}
    x_{i+1,j}^p(t) = \frac{T(T-t)+1}{T(T-i)+1}\tilde x_{i,j}^p + (t-i)\epsilon \bar{Z}_p^{i}, \, t \in [i,i+1].
\end{align*}
Hence, by definition, for $0 \le i \le T-2$,
\begin{equation*}
    \tilde x_{i+1,j}^p = x_{i+1,j}^p(i+1) =  \frac{T(T-i-1)+1}{T(T-i)+1}\tilde x_{i,j}^p + \epsilon \bar{Z}_p^{i}.
\end{equation*}
Utilizing the above recursive relationship, we obtain, for $0 \le i \le T-1$\footnote{In this section, we take by convention that summation $\sum_{k=m}^nc_k=0$ if $m>n$.},
\begin{equation}\label{thm1:eq1}
    \tilde x_{i,j}^p = \frac{T(T-i) +1}{T^2+1}\tilde x_{0,j}^p + \sum_{k=0}^{i-1}\frac{T(T-i)+1}{T(T-k-1)+1}\epsilon\bar{Z}_p^k.
\end{equation}

\paragraph*{Distribution Difference.}
For $0 \le i \le T-1$ and $t \in [i,i+1)$, using the control~\eqref{hat_u_v}, we have
\begin{equation*}
    \dot{x}_p(t) = -\frac{T}{T(T-t)+1}x_p(t) + \frac{T(T-i)+1}{T(T-t)+1}\epsilon \bar{Z}_i^p.
\end{equation*}
Solving the above ODE with the initial condition $x_p(0) = \txinit$, we can get the solution
\begin{equation}\label{IVP_a}
    x_p(t)= \frac{T(T-t)+1}{T^2+1}\txinit + \sum_{k=0}^{i-1}\frac{T(T-t)+1}{T(T-k-1)+1}\epsilon \bar{Z}_p^k + (t-i)\epsilon \bar{Z}_i^p,
\end{equation}
when $0 \le i \le T-1$ and $t \in [i,i+1)$. The above equation also holds when $i = T-1$ and $t = T$.

On the other hand, combining equations \eqref{eq:adaptive_last_iteration} and \eqref{thm1:eq1}, we know that when $0 \le i \le T-1$ and $t \in [i,i+1)$ or $i = T-1$ and $t = T$,
\begin{align*}
    \hat{x}_{T-1,j}^p(t) &= \frac{T(T-t)+1}{T(T-i)+1}\tilde x_{i,j}^p + (t-i)\epsilon Z_{i,j}^p \\ 
    &= \frac{T(T-t)+1}{T^2+1}\tilde x_{0,j}^p + \sum_{k=0}^{i-1}\frac{T(T-t)+1}{T(T-k-1)+1}\epsilon\bar{Z}_p^k + (t-i)\epsilon Z_{i,j}^p.
\end{align*}
The above equation also holds when $i = T-1$ and $t = T$.

Combining the last equation with equation \eqref{IVP_a} and the fact that $\{Z_{i,j}^p\}_{0\le i \le T-1, 1\le j \le N}$, $\tilde{x}_{0,j}^p$ and $\txinit$ are independent normal random variables and 
\begin{equation*}
    \tilde x_{0,j}^p, \txinit \sim \mathcal{N}(0,1), Z_{i,j}^p \sim \mathcal{N}(m,\sigma^2),
\end{equation*}
we know that $\hat{x}_{T-1,1}^p(t)$ and $x_p(t)$ are normal random variables with
\begin{equation*}
\mathbb{E}\hat{x}_{T-1,1}^p(t) = \mathbb{E} x_p(t),    |\Var|\hat{x}_{T-1,1}^{p}(t)|^2 - \Var|x_p(t)|^2| = \sigma^2\epsilon^2(t-i)^2(1-\frac{1}{N}) \le \sigma^2\epsilon^2.
\end{equation*}

\paragraph*{Performance Difference.} Recalling equation \eqref{IVP_a} and plugging \eqref{IVP_a} into \eqref{hat_u_v}, we know that
\begin{align*}
&x_p(T) = \frac{1}{T^2+1}\txinit + \sum_{k=0}^{T-1}\frac{1}{T(T-k-1)+1}\epsilon \bar{Z}_p^k,\\
    &u_p(t) = -\frac{T}{T^2+1}\txinit - \sum_{k=0}^{i-1}\frac{T}{T(T-k-1)+1}\epsilon \bar{Z}_p^k +\epsilon \bar{Z}_i^p, \; 0 \le i \le T-1, \, t \in [i,i+1).
\end{align*}

To compute the difference between $J_p$ and $J_o$, we first notice that
\begin{align*}
    \mathbb{E}J_p - J_o &= \left[\frac{1}{T}\int_{0}^{T}\mathrm{Var}(u_p(t))\rmd t + \mathrm{Var}(x_p(T))\right] + \left[\frac{1}{T}\int_{0}^{T}|\mathbb{E} u_p(t)|^2\rmd t + |\mathbb{E} x_p(T)|^2 - J_o\right] \\
    &:= \mathrm{I}_1 + \mathrm{I}_2.
\end{align*}

By the independence of $\{\bar{Z}_i^p\}_{i=0}^{T-1}$, we know that
\begin{align*}
 \mathrm{I}_1 &= \frac{\sigma^2\epsilon^2}{NT}\sum_{i=0}^{T-1}\sum_{k=0}^{i-1}\frac{T^2}{[T(T-k-1)+1]^2} + \frac{\sigma^2\epsilon^2}{N} + \sum_{k=0}^{T-1}\frac{1}{[T(T-k-1)+1]^2}\frac{\sigma^2\epsilon^2}{N} \\
 &=\frac{\sigma^2\epsilon^2}{N}\left[1 + \sum_{i=0}^{T-1}\sum_{k=0}^{i-1}\frac{T}{[T(T-k-1)+1]^2} + \sum_{k=0}^{T-1}\frac{1}{[T(T-k-1)+1]^2}\right] \\
 &= \frac{\sigma^2\epsilon^2}{N}\left[1 + \sum_{k=0}^{T-2}\sum_{i=k+1}^{T-1}\frac{T}{[T(T-k-1)+1]^2} + \sum_{k=0}^{T-1}\frac{1}{[T(T-k-1)+1]^2}\right] \\
 &= \frac{\sigma^2\epsilon^2}{N}\left[1 + \sum_{k=0}^{T-1}\frac{T(T-k-1)+1}{[T(T-k-1)+1]^2}\right] \\
 &= \frac{\sigma^2\epsilon^2}{N}\left[1 + \sum_{k=0}^{T-1}\frac{1}{Tk+1}\right] \\
 &\le \frac{3\sigma^2\epsilon^2}{N}.
\end{align*}

Meanwhile, noticing that
\begin{align*}
&\mathbb{E}x_p(T) = \frac{1}{T^2+1}\txinit + \sum_{k=0}^{T-1}\frac{1}{T(T-k-1)+1}\epsilon m, \\
&\mathbb{E}u_p(t) = -\frac{T}{T^2+1}\txinit - \sum_{k=0}^{i-1}\frac{T}{T(T-k-1)+1}\epsilon m +\epsilon m,
\end{align*}
it is straightforward to compute that
\begin{equation*}
    \mathrm{I}_2 =\frac{2\epsilon m \txinit}{T^2+1} \mathrm{I}_3 + \epsilon^2 m^2 (\mathrm{I}_4+1),
\end{equation*}
where
\begin{align*}
    \mathrm{I}_3 &=  \sum_{k=0}^{T-1}\frac{1}{T(T-k-1)+1} + \sum_{i=0}^{T-1}\sum_{k=0}^{i-1}\frac{T}{T(T-k-1)+1} - T\\
    &=  \sum_{k=0}^{T-1}\frac{1}{T(T-k-1)+1}+ \sum_{k=0}^{T-1}\sum_{i=k+1}^{T-1}\frac{T}{T(T-k-1)+1} - T \\
    &= \sum_{k=0}^{T-1}\frac{T(T-k-1)+1}{T(T-k-1)+1} - T = 0,
\end{align*}
and
\begin{align*}
    \mathrm{I}_4 &=  (\sum_{k=0}^{T-1}\frac{1}{T(T-k-1)+1})^2 + T\sum_{i=0}^{T-1}(\sum_{k=0}^{i-1}\frac{1}{T(T-k-1)+1})^2 - \sum_{i=0}^{T-1}\sum_{k=0}^{i-1}\frac{2}{T(T-k-1)+1}\\
    &=  \sum_{k=0}^{T-1}\frac{1+T(T-k-1)}{[1+T(T-k-1)]^2} +2\sum_{i=0}^{T-1}\sum_{k=0}^{i-1}\frac{T(T-i-1)+1}{[T(T-k-1)+1][T(T-i-1)+1]} \\
    &\qquad\qquad\qquad\qquad\qquad\qquad-2\sum_{i=0}^{T-1}\sum_{k=0}^{i-1}\frac{1}{T(T-k-1)+1} \\
    &= \sum_{k=0}^{T-1}\frac{1}{1+T(T-k-1)} = \sum_{k=0}^{T-1}\frac{1}{Tk+1} \le 2.
\end{align*}
Therefore,
\begin{equation*}
    \mathbb{E} J_p - J_o = \mathrm{I}_1 +\frac{2\epsilon m \txinit}{T^2+1} \mathrm{I}_3 + \epsilon^2 m^2 (\mathrm{I}_4+1) \le 3(m^2+\frac{\sigma^2}{N})\epsilon^2.
\end{equation*}

\subsection{\dagger}
With the same approach of computing $u_v$ in \eqref{eq:vanilla_control}, we have
\begin{equation*}
    u_0(t,x) = -\frac{T}{T(T-t)+1}x + \frac{T^2+1}{T(T-t)+1}\epsilon \bar{Z}_0^d,
\end{equation*}
where
\begin{equation*}
    \bar{Z}_0^d = \frac{1}{N}\sum_{j=1}^{N}Z_{0,j}^d.
\end{equation*}
Recalling the definition of $x_{0,j}^d$ in equation \eqref{ivp_DAGGER}, we have
\begin{equation*}
    x_{0,j}^d(t) = \frac{T(T-t)+1}{T^2+1}\tilde{x}_{0,j}^d + \epsilon t \bar{Z}_0^d.
\end{equation*}
Hence, for $0\le i \le T-1$, we have
\begin{equation*}
    \tilde{x}_{i,j}^d = x_{0,j}^d(i) = \frac{T(T-i)+1}{T^2+1}\tilde{x}_{0,j}^d + \epsilon i \bar{Z}_0^d.
\end{equation*}
Plugging the last equation into equation \eqref{DAGGER_data}, we have that for $t \in [i,T]$
\begin{align*}
    &\hat{u}_{i,j}^d(t) = - \frac{T}{T^2+1}\tilde{x}_{0,j}^d- \epsilon \frac{Ti}{F(i)}\bar{Z}_0^d + \epsilon Z_{i,j}^d,\\ &\hat{x}_{i,j}^d(t) = \frac{T(T-t)+1}{T^2+1}\tilde{x}_{0,j}^d + \epsilon\frac{iF(t)}{F(i)}\bar{Z}_0^d + (t-i)\epsilon Z_{i,j}^d,
\end{align*}
where $F(t) = T(T-t)+1$.
Therefore,
\begin{equation*}
    \hat{u}_{i,j}^d(t) = -\frac{T}{F(t)}\hat{x}_{i,j}^d(t) + \frac{F(i)}{F(t)}\epsilon Z_{i,j}^d.
\end{equation*}
We can then compute the least squares problem \eqref{least_square_DAGGER} to obtain that for $0 \le i \le T-1$ and $t \in [i,i+1)$, we have
\begin{equation*}
    u_d(t,x) = -\frac{T}{F(t)}x + \frac{1}{i+1}\sum_{k=0}^{i}\frac{F(k)}{F(t)}\epsilon \bar{Z}_k^d,
\end{equation*}
where
\begin{equation*}
    \bar{Z}_i^d = \frac{1}{N}\sum_{j=1}^{N} Z_{i,j}^d,
\end{equation*}
for $0 \le i \le T-1$. Hence, we have that when $0 \le i \le T-1$ and $t \in [i,i+1)$
\begin{align*}
    x_d(t) = \frac{T(T-t)+1}{T^2+1}\txinit + \sum_{k=0}^{i-1}\frac{F(t)}{F(k+1)F(k)}\frac{1}{k+1}\sum_{l = 0}^k\epsilon F(l)\bar{Z}_l^d +  \frac{t-i}{(i+1)F(i)}\sum_{k=0}^iF(k)\epsilon \bar{Z}_k^d.
\end{align*}
Therefore, when $0 \le i \le T-1$ and $t \in [i,i+1)$,
\begin{equation*}
    u_d(t) = -\frac{T}{T^2+1}\txinit - \sum_{k=0}^{i-1}\frac{T}{(k+1)F(k)F(k+1)}\sum_{l=0}^k \epsilon F(l)\bar{Z}_l^d +  \frac{1}{(i+1)F(i)}\sum_{k=0}^i F(k) \epsilon \bar{Z}_k^d.
\end{equation*}

\paragraph*{Distribution Difference.}
We first give a lower bound of $\Var(\hat{x}^d(t))$. For any $t \in [i,i+1)$, $\hat{x}^d(t)$ can be viewed as the random variable which uniformly samples an element from $\{\hat{x}_{k,j}(t)\}_{0 \le k \le i, 1 \le j \le N}$. Noticing that $\Var(X) = \Var(\mathbb{E}(X|Y)) + \mathbb{E}\Var(X|Y)$, we know that for any $t \in [i,i+1)$
\begin{align*}
    \Var(\hat{x}^d(t)) \ge \frac{1}{N(i+1)}\sum_{k=0}^i\sum_{j=1}^N \Var(\hat{x}_{k,j}^d(t)).
\end{align*}
Then, through the independence of $\tilde{x}_{0,j}^d$, $Z_{0,j}^d$ and $Z_{i,j}^d$, we have that
\begin{align*}
     \Var(\hat{x}^d(t)) &\ge \left[\frac{T(T-t)+1}{T^2+1}\right]^2 + \frac{1}{i+1}\sum_{k=0}^i (t-k)^2\epsilon^2 \sigma^2 \\
     &= \left[\frac{T(T-t)+1}{T^2+1}\right]^2 + \epsilon^2\sigma^2\left[t^2 - ti + \frac{i(2i+1)}{6}\right].
\end{align*}
Note that
\begin{equation*}
    t^2 - ti + \frac{i(2i+1)}{6} - \frac{t^2}{3} = t(t-i) + \frac{i}{6} - \frac{1}{3}(t+i)(t-i) \ge (2/3 t - 1/3 i)(t-i) \ge 0.
\end{equation*}
Hence,
\begin{equation}\label{thm:eq11}
    \Var(\hat{x}^d(t)) \ge\left[\frac{T(T-t)+1}{T^2+1}\right]^2 + \frac{t^2\epsilon^2\sigma^2}{3}.
\end{equation}

Next, we give an upper bound of $\Var(x_d(t))$. In this part, without loss of generality, we assume $m = 0$, as the value of $m$ does not influence of the value of $\Var(x_d(t))$. Let 
\begin{align*}
    x'_d(t) &= x_d(t) - \frac{T(T-t)+1}{T^2+1}\txinit \\
    &= \sum_{k=0}^{i-1}\frac{F(t)}{F(k+1)F(k)}\frac{1}{k+1}\sum_{l = 0}^k\epsilon F(l)\bar{Z}_l^d +  \frac{t-i}{(i+1)F(i)}\sum_{k=0}^iF(k)\epsilon \bar{Z}_k^d.
\end{align*}
Then, 
\begin{equation}\label{thm:eq12}
    \Var(x_d(t))=\left[\frac{T(T-t)+1}{T^2+1}\right]^2 + \Var(x'_d(t)).
\end{equation}
and for any $t \in [i,i+1)$, $\dot{x}'_d(t) = u_d(t,x'_d(t))$. Define $V(t) = \Var(x'_d(t)) = \mathbb{E} |x'_d(t)|^2$. Then, we have that
\begin{align*}
    \dot{V}(t) = 2 \mathbb{E} x'_d(t) u_d(t,x'_d(t))
    &= -\frac{2T}{F(t)}V(t) + 2\frac{\epsilon^2\sigma^2}{N}\mathrm{I}_5, 
\end{align*}
where
\begin{equation*}
    \mathrm{I}_5 = \sum_{k=0}^{i-1}\frac{1}{F(k+1)F(k)(k+1)(i+1)}\sum_{l = 0}^k F^2(l) +  \frac{t-i}{(i+1)^2F(i) F(t)}\sum_{k=0}^iF^2(k).
\end{equation*}
Noticing that $F(l) \le F(0)$, we have
\begin{align*}
    I_5 &\le \frac{F^2(0)}{i+1}\left[\sum_{k=0}^{i-1}\frac{1}{F(k+1)F(k)} + \frac{(t-i)}{F(i)F(t)}\right] \\
    &= \frac{F^2(0)}{T(i+1)}\left[\sum_{k=0}^{i-1}(\frac{1}{F(k+1)} - \frac{1}{F(k)}) + \frac{1}{F(t)} - \frac{1}{F(0)} \right] \\
    &=\frac{F^2(0)}{T(i+1)}\left[\frac{1}{F(t)} - \frac{1}{F(0)}\right] = \frac{F^2(0)}{T(i+1)}\frac{Tt}{F(0)F(t)} \le \frac{F(0)}{F(t)}.
\end{align*}
Therefore, 
\begin{equation*}
    \dot{V}(t) \le -\frac{2T}{F(t)}V(t) + \frac{2\epsilon^2\sigma^2}{N}\frac{F(0)}{F(t)}.
\end{equation*}
Let $W(t) = 2\frac{\epsilon^2\sigma^2}{N} t - V(t)$, then $W(0) = 0$ and
\begin{align*}
    \dot{W}(t) + \frac{2T}{F(t)} W(t) &\ge \frac{\epsilon^2\sigma^2}{N}[2 + \frac{4Tt}{F(t)}] -\dot{V}(t) - \frac{2T}{F(t)} \\
    &\ge \frac{\epsilon^2\sigma^2}{N}[2+\frac{4Tt}{F(t)} - 2\frac{F(0)}{F(t)}] \\
    &= \frac{\epsilon^2\sigma^2}{N}[2+\frac{4Tt}{F(t)} - 2\frac{2F(0)}{F(t)}] = \frac{\epsilon^2\sigma^2}{N}\frac{2Tt}{F(t)}\ge 0.
\end{align*}
Hence, $W(t) \ge 0$, which means that
\begin{equation*}
    \Var(x'_d(t)) = V(t) \le 2\frac{\epsilon^2\sigma^2}{N}t
\end{equation*}
Combining the last inequality with \eqref{thm:eq11} and \eqref{thm:eq12}, we have
\begin{equation*}
    \Var(\hat{x}^d(t)) - \Var(x_d(t)) \ge \epsilon^2\sigma^2(\frac{t^2}{3}-\frac{2t}{N}).
\end{equation*}

\paragraph*{Performance Difference. } Define
\begin{equation*}
    e_i = - \sum_{k=0}^{i-1}\frac{T}{(k+1)F(k)F(k+1)}\sum_{l=0}^k \epsilon F(l)\bar{Z}_l^d +  \frac{1}{(i+1)F(i)}\sum_{k=0}^i F(k) \epsilon \bar{Z}_k^d,\quad 0 \le i \le T-1.
\end{equation*}
We have
\begin{align*}
     J_d &= \frac{1}{T}\sum_{i=0}^{T-1}|-\frac{T}{T^2+1}\txinit +e_i|^2 + |\txinit - \frac{T^2}{T^2+1}\txinit + \sum_{i=0}^{T-1}e_i|^2 \\
     &= \frac{T^2|\xinit|^2}{(T^2+1)^2} - \sum_{i=0}^{T-1}\frac{2\txinit e_i}{T^2+1} + \frac{1}{T}\sum_{i=0}^{T-1}|e_i|^2 + \frac{|\txinit|^2}{(T^2+1)^2} + \sum_{i=0}^{T-1}\frac{2e_i\txinit}{T^2+1} + |\sum_{i=0}^{T-1}e_i|^2 \\
     &= \frac{|\txinit|^2}{T^2+1} + \frac{1}{T}\sum_{i=0}^{T-1}|e_i|^2 + |\sum_{i=0}^{T-1}e_i|^2.
\end{align*}
Therefore
\begin{equation}\label{thm2_eq1}
    \mathbb{E} J_d - J_o \ge \mathbb{E} |\sum_{i=0}^{T-1}e_i|^2 = (\mathbb{E}\sum_{i=0}^{T-1}e_i)^2 + \mathrm{Var}(\sum_{i=0}^{T-1}e_i).
\end{equation}
We can then compute that
\begin{align*}
    \sum_{i=0}^{T-1}e_i &= -\sum_{i=0}^{T-1}\sum_{k=0}^{i-1}\sum_{l=0}^k\frac{\epsilon TF(l)\bar{Z}_l^d}{(k+1)F(k)F(k+1)}  + \sum_{i=0}^{T-1}\sum_{k=0}^i \frac{F(k) \epsilon \bar{Z}_k^d}{(i+1)F(i)}\\
    &=\sum_{i=0}^{T-1}\sum_{k=0}^i \frac{F(k) \epsilon \bar{Z}_k^d}{(i+1)F(i)} - \sum_{k=0}^{T-1}\sum_{l=0}^k\sum_{i=k+1}^{T-1}\frac{\epsilon TF(l)\bar{Z}_l^d}{(k+1)F(k)F(k+1)} \\
     &=\sum_{i=0}^{T-1}\sum_{k=0}^i \frac{F(k) \epsilon \bar{Z}_k^d}{(i+1)F(i)} - \sum_{i=0}^{T-1}\sum_{k=0}^i\frac{\epsilon T(T-i-1)F(k)\bar{Z}_k^d}{(i+1)F(i)F(i+1)} \\
     &= \sum_{i=0}^{T-1}\sum_{k=0}^i\frac{F(k)\epsilon \bar{Z}_k^d}{(i+1)F(i)F(i+1)}.
\end{align*}
Therefore,
\begin{equation}\label{thm2_eq2}
\begin{aligned}
        \mathbb{E}\sum_{i=0}^{T-1}e_i &= \epsilon m \sum_{i=0}^{T-1}\sum_{k=0}^i\frac{F(k)}{(i+1)F(i)F(i+1)}\\
    &= \epsilon m\sum_{i=0}^{T-1}\frac{(T^2+1)(i+1)-Ti(i+1)/2}{(i+1)F(i)F(i+1)} \\
    &= \epsilon m\sum_{i=0}^{T-1}\frac{T^2+1 - Ti/2}{[T(T-i)+1][T(T-i-1)+1]} \\
    &\ge \epsilon m\frac{T^2+T+2}{2T}\sum_{i=0}^{T-1}[\frac{1}{T(T-i-1)+1} - \frac{1}{T(T-i)+1}]\\
    &= \epsilon m\frac{T^2+T+2}{2T}(1-\frac{1}{T^2+1}) \ge \frac{\epsilon m T}{2}. 
\end{aligned}
\end{equation}
  On the other hand, noticing that
\begin{equation*}
    \sum_{i=0}^{T-1}e_i =  \sum_{k=0}^{T-1}\sum_{i=k}^{T-1}\frac{F(k)\epsilon \bar{Z}_k^d}{(i+1)F(i)F(i+1)},
\end{equation*}
we have
\begin{equation}\label{thm2_eq3}
\begin{aligned}
        \mathrm{Var}(\sum_{i=0}^{T-1}e_i) &= \frac{\epsilon^2 \sigma^2}{N}\sum_{k=0}^{T-1}F^2(k)(\sum_{i=k}^{T-1}\frac{1}{(i+1)F(i)F(i+1)})^2\\
    &\ge \frac{\epsilon^2\sigma^2}{NT^4}\sum_{k=0}^{T-1}[T(T-k)+1]^2(\sum_{i=k}^{T-1}\frac{1}{T(T-i-1)+1} -\frac{1}{T(T-i)+1})^2 \\
    &= \frac{\epsilon^2\sigma^2}{NT^4}\sum_{k=0}^{T-1}[T(T-k)+1]^2[1- \frac{1}{T(T-k)+1}]^2 \\
    &= \frac{\epsilon^2\sigma^2}{NT^2}\sum_{k=0}^{T-1}(T-K)^2 = \frac{\epsilon^2\sigma^2}{NT^2}\frac{T(T+1)(2T+1)}{6} \ge \frac{\epsilon^2\sigma^2T}{3N}.
\end{aligned}
\end{equation}
Combining equations \eqref{thm2_eq1}, \eqref{thm2_eq2} and \eqref{thm2_eq3}, we conclude our result.

\section{Full Dynamics of Quadrotor}
\label{appendix:dyn}

In this section, we introduce the full dynamics of quadrotor~\citep{bouabdallah2004design,madani2006control,mahony2012multirotor} that are considered in Section \ref{sec:prob_quadrotor}. %
The state variable of a quadrotor is $\bm x=(\bm p\transpose,\bm v_b\transpose,\bm \eta\transpose, \bm w_b\transpose)\transpose\in \mathbb{R}^{12}$ where $\bm p=(x,y,z)\in\mathbb{R}^3$ is the position of quadrotor in Earth-fixed coordinates, $\bm v_b\in\mathbb{R}^3$ is the velocity in body-fixed coordinates, $\bm \eta=(\phi, \theta, \psi )\in\mathbb{R}^3$ (roll, pitch, yaw) is the attitude in terms of Euler angles in Earth-fixed coordinates, and $\bm w_b\in\mathbb{R}^3$ is the angular velocity in body-fixed coordinates. Control $\bm u=( s,\tau_x,\tau_y,\tau_z)\transpose\in \mathbb{R}^4$ is composed of total thrust $s$ and body torques $(\tau_x,\tau_y,\tau_z)$ from the four rotors. Then we can model the quadrotor's dynamics as
\begin{equation}\label{eq:uav_old1}
  \begin{cases}
    \dot{\bm{p}}= \bm{R}\transpose(\bm{\eta})\bm{v}_b,\\
    \dot{\bm{v}}_b= - \bm{w}_b\times \bm{v}_b - \bm{R}(\bm{\eta})\bm{g} + \frac{1}{m}A\bm u,\\
    \dot{\bm{\eta}} = \bm{K}(\bm{\eta})\bm{w}_b, \\
    \dot{\bm{w}}_b = - \bm{J}^{-1}\bm{w}_b\times \bm{J}\bm{w}_b + \bm{J}^{-1}B \bm u,
  \end{cases}
\end{equation}
with matrix $A$ and $B$  defined as
\begin{equation*}
  A
  =
  \begin{bmatrix}
    0 & 0 & 0 & 0 \\
    0 & 0 & 0 & 0 \\
    1 & 0 & 0 & 0
  \end{bmatrix},
  \quad\quad
  B
  =
  \begin{bmatrix}
    0 & 1 & 0 & 0 \\
    0 & 0 & 1 & 0 \\
    0 & 0 & 0 & 1
  \end{bmatrix}.
\end{equation*}
The constant mass $m$ and inertia matrix $\bm J=\text{diag}(J_{x}, J_{y}, J_{z})$ are the parameters of the quadrotor, where $J_{x}, J_{y}$, and $J_{z}$ are the moments of inertia of the quadrotor in the $x$-axis, $y$-axis, and $z$-axis, respectively. We set $m=2kg $ and $J_x=J_y=\frac{1}{2}J_z=1.2416 kg\cdot m^2$ which are the same system parameters as in \citep{madani2006control}. The constants  $\bm{g}= (0,0, g)\transpose$ denote the gravity vector where $g=9.81 m/s^2$ is the acceleration of gravity on Earth. 
The direction cosine matrix $\bm{R}(\bm{\eta})\in SO(3)$ represents the transformation from the Earth-fixed coordinates to the body-fixed coordinates: 
\begin{equation*}
  \bm{R}(\bm{\eta}) =
  \begin{bmatrix}
    \cos{\theta}\cos{\psi} & \cos{\theta}\sin{\psi} & -\sin{\theta}\\
    \sin{\theta}\cos{\psi}\sin{\phi}-\sin{\psi}\cos{\phi} &
    \sin{\theta}\sin{\psi}\sin{\phi}+\cos{\psi}\cos{\phi} &
    \cos{\theta}\sin{\phi} \\
    \sin{\theta}\cos{\psi}\cos{\phi}+\sin{\psi}\sin{\phi} &
    \sin{\theta}\sin{\psi}\cos{\phi}-\cos{\psi}\sin{\phi} &
    \cos{\theta}\cos{\phi}
  \end{bmatrix},
\end{equation*}
and the attitude kinematic matrix $\bm{K}(\bm{\eta})$ relates the time derivative of the attitude representation with the associated angular rate:
\begin{equation*}
  \bm{K}(\bm{\eta}) =
  \begin{bmatrix}
    1 & \sin{\phi}\tan{\theta} & \cos{\phi}\tan{\theta}\\
    0 & \cos{\phi} & -\sin{\phi} \\
    0 & \sin{\phi}\sec{\theta} & \cos{\phi}\sec{\theta}
  \end{bmatrix},
\end{equation*}

Note that in practice the quadrotor is directly controlled by the individual rotor thrusts $\bm{F}=(F_1, F_2, F_3, F_4)\transpose$, and we have the relation $\bm{u}=E\bm{F}$ with
\begin{equation*}
  E
  =
  \begin{bmatrix}
    1 & 1 & 1 & 1  \\
    0 & l & 0 & -l \\
    -l & 0 & l & 0 \\
    c & -c & c & -c
  \end{bmatrix},
\end{equation*}
where $l$ is the distance from the rotor to the UAV's center of gravity and $c$ is a constant that relates the rotor angular momentum to the rotor thrust (normal force). 
So once we obtain the optimal control $\bm{u}^*$, we are able to get the optimal $\bm{F}^*$ immediately by the relation $\bm{F}^*=E^{-1}\bm{u}^*$.

\section{PMP and Space Marching Method}
\label{appendix:PMP}
In this section we introduce the open-loop optimal problem solver used for solving the optimal landing problem of a quadrotor. The solver is based on Pontryagin's minimum principle (PMP) \citep{pontryagin1987mathematical} and space-marching method \citep{zang2022machine}.
The  optimal landing problem is defined as
\begin{equation*}
\begin{aligned}
 \min_{\bm{x},\bm{u}}   &\int_{0}^{T}L(\bm{x}(\tau),\bm{u}(\tau))d\tau+M(\bm{x}(T)), \\ 
 \text{s.t.}& 
 \begin{cases}
 \dot{\bm{x}}(t)=f(\bm{x}(t),\bm{u}(t)), t\in[0,T], \\
 \bm{x}(0)=\bm{x}_0, \\
 \end{cases}
\end{aligned}
 \label{eq:quad}
\end{equation*}
where $\bm{x}(t):[0,T]\rightarrow \mathbb{R}^{12}$ and $\bm{u}(t):[0,T]\rightarrow \mathbb{R}^{4}$ denote the state trajectory and control trajectory, respectively, and $f$ is the full dynamics of quadrotor introduced in Appendix \ref{appendix:dyn}.
By PMP, problem \eqref{eq:quad} can be solved through solving  a two-point boundary value problem (TPBVP). Introduce costate variable  $\bm{\lambda}\in \mathbb{R}^{12}$ and Hamiltonian
\begin{equation*}
H (\bm x, \bm \lambda, \bm u) = L(\bm x, \bm u) + \bm \lambda \cdot f (\bm x, \bm u).
\end{equation*}
 The TPBVP is defined as
\begin{equation}
\label{TP1}
\left\{\begin{array}{l}
\dot{\boldsymbol{x}}(t)=\partial_{\boldsymbol{\lambda}}^{T} H\left(\boldsymbol{x}(t), \boldsymbol{\lambda}(t), \boldsymbol{u}^{*}(t)\right), \\
\dot{\boldsymbol{\lambda}}(t)=-\partial_{\boldsymbol{x}}^{T} H\left(\boldsymbol{x}(t), \boldsymbol{\lambda}(t), \boldsymbol{u}^{*}(t)\right).
\end{array}\right.
\end{equation}
We have the boundary conditions:
\begin{equation}
\label{TP2}
    \left\{\begin{array}{l}
\boldsymbol{x}(0)=\boldsymbol{x}_{0}, \\
\boldsymbol{\lambda}\left(T\right)=
\nabla M\left(\boldsymbol{x}\left(T\right)\right),
\end{array}\right.
\end{equation}
and the optimal control $\boldsymbol{u}^{*}(t)$ should  minimize Hamiltonian at each $t$ :
\begin{equation}
\label{TP3}
  \boldsymbol{u}^{*}(t)=\underset{\boldsymbol{u} }{\arg \min }\ \  H(\boldsymbol{x}(t), \boldsymbol{\lambda}(t), \boldsymbol{u}).
\end{equation}

We use \textit{solve\_bvp} function of \textit{scipy} \citep{kierzenka2001bvp} to solve TPBVP \eqref{TP1}-\eqref{TP3} and set \textit{tolerance} to $10^{-5}$, \textit{max\_nodes} to $5000$. We note that when the initial state $\bm x_0$ is far from the target state $\bm x_T$, solving the TPBVP directly often fails.
Thus we use the space-marching method proposed in \citet{zang2022machine}. We uniformly select $K$ points in the line segment from $\bm x_T$ to $\bm x_0$,
and denote them as $\{\bm x^1_0, \bm x^2_0, \cdots,\bm  x^K_0\}$  according to their increasing distances to $\bm x_T$ ($x^K_0 = \bm x_0$). These $K$ TPBVPs will be solved in order and at every step we use the previous solution as the initial guess to the current problem.  \revisionjmla{With this strategy, every open-loop trajectory in our experiments converged. In later runs, we warm-start the solver with the IVP trajectory generated by our trained neural controller, which almost eliminates the need for space-marching.}

\section{Details on the 7-DoF Manipulator}\label{app:dynamics_manipulator}

In this section, we introduce the dynamics for the 7-DoF torque-controlled manipulator.
Recall that the dynamics of the manipulator is, 
\[
    \dot{\bm x} = \bm f(\bm x, \bm u) = (\bm v, \bm a(\bm x, \bm u))
,\] 
where $\bm u\in\bR^7$ is the control torque, $\bm x = (\bm q, \bm v) \in\bR^{14}$, $\bm q\in\bR^7$ is the joint angles, $\bm  v=\dot{\bm q}\in\bR^{7}$ is the joint velocities, $\ddot{\bm q} = \bm a(\bm x, \bm u) \in \bR^7 $ is the acceleration of joint angles.
To close the equation, we write down the inverse dynamics of the manipulator,
\begin{equation*}
    M(\bm q)\bm a +  C(\bm q, \dot{\bm q})\dot{\bm q} + \bm g(\bm q)= \bm u,
\end{equation*}
where one can compute the acceleration in terms of $\bm x$ and $\bm u$.
Here $M(\bm q)$ is the generalized inertia matrix, $C(\bm q,\dot{\bm q})\dot{\bm q}$ represents the \textit{centrifugal} forces and \textit{Coriolis} forces, and $\bm g(\bm q)$ is the generalized gravity.

\section{The QRNet}\label{sec:QRnet}
In this section, we introduce the structure of QRNet. QRNet utilize the linear quadratic regulator (LQR) controller at an equilibrium and thus improve stability around the equilibrium.
Suppose we have the LQR controller, $\bm u^{\LQR}$, for the problem with linearized dynamics and quadratized costs at $(\bm x_1, \bm u_1)$, the QRNet can be formulated as
\begin{equation*}
    \bm u^{\QR}(t,\bm x) = \sigma (\bm u^{\LQR}(t, \bm x) + \hat{\bm {u}}(t, \bm x;\theta) - \hat{\bm {u}}(T,\bm x_1;\theta)), 
\end{equation*}
where $\hat{\bm {u}}(t,\bm x;\theta) $ is any neural network with trainable parameters $\theta$, and  $\sigma$ is a saturating function that satisfies $\sigma(\bm u_1) = \bm u_1, \sigma_{\bm u}(\bm u_1)=I$, where $I$ is the identity matrix.
The $\sigma$ used in this example is defined coordinate-wisely as
\begin{equation*}
    \sigma (u) = u_\tmin + \frac{u_\tmax - u_\tmin}{1+c_1 \exp[-c_2(u-u_1)]},
\end{equation*}
where
$ c_1=(u_\tmax - u_1) / (u_1-u_\tmin), c_2=(u_\tmax - u_\tmin) / [(u_\tmax -u_1)(u_1-u_\tmin)]$ with
$u_\tmin, u_\tmax$ being minimum and maximum values for $u$.
Here $u, u_\tmin$ and $u_\tmax$ are the corresponding values at each coordinate of $\bm u, \bm u_\tmin, \bm u_\tmax$, respectively.
In the first experiment evaluated in Section \ref{sec:reaching_problem}, we set $u_\tmin = -150$ and $u_\tmax = 15$. In the more challenging scenario in Section \ref{sec:reaching_compare}, we adjust the bounds to $u_\tmin = -2000$ and $u_\tmax = 2000$ to prevent saturation, as \dagger explores a wider range of states, leading to larger control torques.

To get the $\bm u^{\LQR}$, we expand the dynamics linearly as
\begin{equation*}
    \bm f(\bm x, \bm u) \approx \bm f_{\bm x}(\bm x_1, \bm u_1) (\bm x - \bm x_1) + \bm f_{\bm u}(\bm x_1, \bm u_1) (\bm u - \bm u_1),
\end{equation*}
and the term related to acceleration in the running cost quadratically as
\begin{align*}
    &\bm a(\bm x,\bm u)\transpose Q_{\bm a} \bm a(\bm x,\bm u)
    \approx ~  \mathcal{L}_{\bm a}(\bm x,\bm u)\transpose Q_{\bm a} \mathcal{L}_{\bm a}(\bm x, \bm u) \\
    =  &(\bm x-\bm x_1)\transpose \bm a_{\bm x} \transpose Q_{\bm a}\bm a_{\bm x}(\bm x-\bm x_1)+ (\bm u-\bm u_1)\transpose \bm a_{\bm u}\transpose Q_{\bm a} \bm a_{\bm u}(\bm u-\bm u_1) 
    + 2(\bm x-\bm x_1)\transpose \bm a_{\bm x}\transpose Q_{\bm a} \bm a_{\bm u}(\bm u-\bm u_1),
\end{align*}
where $\mathcal{L}_{\bm a} = \bm a_{\bm x}(\bm x_1, \bm u_1)(\bm x-\bm x_1) + \bm a_{\bm u}(\bm x_1, \bm u_1)(\bm u-\bm u_1)$, and we exploit $\bm a(\bm x_1, \bm u_1) =\bm 0 $ and $\bm f(\bm x_1, \bm u_1) = \bm 0 $. 
The derivatives boil down to $\bm a_{\bm x}$ and $\bm a_{\bm u} $ which can be analytically computed in the Pinocchio library \citep{pinocchioweb,carpentier2019pinocchio,carpentier2018analytical}.
In the experiment, we solve the LQR by the implementation in the Drake library \citep{drake}.

\bibliography{ref}
\bibliographystyle{unsrt}

\end{document}